\documentclass[a4paper, 10pt]{article}
\usepackage{amssymb, amsthm, amsmath}
\numberwithin{equation}{section}

\newcommand{\mcp}{\mathcal{P}}
\newcommand{\R}{\mathbb{R}}
\newcommand{\C}{\mathbb{C}}
\newcommand{\Z}{\mathbb{Z}}
\newcommand{\N}{\mathbb{N}}

\newcommand{\Q}{\mathbb{Q}}
\newcommand{\x}{\ma{x}}
\newcommand{\cQ}{\overline{\mathbb{Q}}}
\newcommand{\Grass}{\mathbb{G}}
\newcommand{\e}{\emph}
\newcommand{\rom}{\mathrm}
\newcommand{\bfP}{\mathbb{P}}
\newcommand{\A}{\mathbb{A}}

\newcommand{\ma}{\mathbf}
\newcommand{\ben}{\begin{enumerate}}
\newcommand{\een}{\end{enumerate}}
\newcommand{\eit}{\begin{itemize}}
\newcommand{\beq}{\begin{equation}}
\newcommand{\eeq}{\end{equation}}

\newcommand{\ve}{\varepsilon}

\newcommand{\mcal}{\mathcal}

\newcommand{\lab}{\label}

\newcommand{\al}{\alpha}
\newcommand{\D}{\Delta}
\newcommand{\del}{\delta}

\newcommand{\be}{\beta}

\newcommand{\la}{\lambda}
\newcommand{\sfl}{\mathsf{\Lambda}}
\newcommand{\sfm}{\mathsf{M}}
\newcommand{\sfn}{\mathsf{N}}
\newcommand{\sfg}{\mathsf{\Gamma}}

\newtheorem{thm}{Theorem}
\newtheorem{lem}{Lemma}

\newtheorem{cor}{Corollary}
\newtheorem*{con}{Conjecture}
\newcommand{\hcf}{\rom{h.c.f.}}
\renewcommand{\mod}{\hspace{-1mm}\pmod}
\newcommand{\colt}[2]{\genfrac{}{}{0pt}{1}{#1}{#2}}

\newcommand{\w}{\ma{w}}
\newcommand{\y}{\ma{y}}

\begin{document}

\title{Counting Rational Points on Hypersurfaces}

\author{T.D. Browning$^1$ \\ D.R. Heath-Brown$^2$\\
\small\emph{Mathematical Institute,
24--29 St. Giles',
Oxford OX1 3LB}\\
\small{$^1$browning@maths.ox.ac.uk}, \small{$^2$rhb@maths.ox.ac.uk}}
\date{}

\maketitle

\begin{abstract}
For any $n \geq 2$, let $F \in \mathbb{Z}[x_1,\ldots,x_n]$ be a form of degree
$d\geq 2$, which produces a geometrically irreducible hypersurface
in $\mathbb{P}^{n-1}$. This paper is concerned with the number
$N(F;B)$ of rational points on $F=0$ which have height at most $B$.
For any $\varepsilon>0$ we establish the estimate 
$$
N(F;B)=O(B^{n-2+\varepsilon}),
$$
whenever either $n\leq 5$ or the hypersurface is not a union of lines.
Here the implied constant depends at most upon $d, n$ and $\varepsilon$.  
\end{abstract}

\section{Introduction}\lab{intro}

For any $n \geq 2$, let $F \in \Z[x_1,\ldots,x_n]$ be a non-zero form of degree
$d$, which produces an algebraic hypersurface $X$ in $\bfP^{n-1}$.  
Our basic interest is with the distribution of rational points on such
hypersurfaces. 
It will be convenient to write $Z^n$ for the 
set of primitive vectors $\x=(x_1,\ldots,x_n) \in\Z^n$, 
where $\x$ is said to be primitive if  $\hcf(x_1,\ldots,x_n)=1.$ 
With this notation, we seek  to understand the asymptotic behaviour
of the quantity 
$$
N(F;B)=\#\{\x \in Z^n: F(\x)=0, ~\mbox{$\max_{i} |x_i| \leq B$}\},
$$
as $B\rightarrow \infty$.  
More precisely, this paper is motivated by the following conjecture of
the second author
 \cite[Conjecture $2$]{annal}.

\begin{con}  Let $\ve>0$.
Suppose that $F$ is  irreducible
and  that $d \geq 2$.
Then we have
$$
N(F;B)=O(B^{n-2+\ve}).
$$
\end{con}

Here, and throughout this paper, the implied constant may
depend at most upon $d, n$ and the choice of $\ve$.  Any further 
dependences will be
explicitly indicated. In the statement of the above conjecture, and
elsewhere, we shall always take irreducibility of
a form to mean absolute irreducibility.
One might also ask about bounds of the shape $N(F;B)=O_F(B^{n-2+\ve})$,
as in \cite[Conjecture $1$]{annal}, where the implied constant is
allowed to depend on the coefficients of $F$. However it transpires
that estimates which are uniform in forms of fixed degree, and with a
fixed number of variables, are much more useful in applications, and
that in most cases results with weaker uniformity appear to be no 
easier to prove.

We take a moment to discuss the available evidence for the conjecture.
Firstly it should be clear that the bound's exponent is in general essentially
as sharp as can be hoped for.
Indeed, whenever the hypersurface $X$ contains a linear
space defined over $\Q$ and having dimension $n-3$, then we
automatically have $N(F;B)\gg_F B^{n-2}$.
The second author has already established the
conjecture in the case of quadrics \cite[Theorem 2]{annal}, in
addition to the cases $n=3$ and $n=4$ for any degree 
\cite[Theorems 3 and 9]{annal}.  
More recently, Broberg and Salberger \cite{broberg} have examined the
case $n=5$ of threefolds. They succeed in establishing the conjecture
as soon as $d \geq 4$.
The best result available in most other cases is the bound 
\beq\lab{pila}
N(F;B)=O(B^{n-2+1/d+\ve}),
\eeq
due to Pila \cite{pila}.  
In the case of cubic threefolds, Broberg and
Salberger \cite{broberg} have improved upon the exponent $3+1/3$ appearing in
(\ref{pila}), allowing one to replace it by $3+1/18$.

Before proceeding further we record the following rather general 
``trivial'' 
bound, in which $[\ma{x}] \in \bfP^{N-1}$ denotes the projective point 
corresponding to a vector $\ma{x}\in \C^{N}$ for any $N \geq 2$. 
This will be established in \S \ref{pre}.

\begin{thm}\label{trivthm}
Let $Y\subseteq\bfP^{N-1}$ be an irreducible variety defined over 
$\cQ$, of dimension $m$ and degree $D$.  Then 
\[\#\{\x\in\Z^N: [\x] \in Y, ~\mbox{$\max_{i} |x_i| \leq B$}\}\ll_{D,N} 
B^{m+1}.\]
\end{thm}

Here, and elsewhere, when we say that a variety is irreducible, we 
shall mean that it is geometrically irreducible. 
We may conclude from Theorem \ref{trivthm} that 
points which lie on a subvariety of $X$ of degree  
$O(1)$, and codimension 
1 or more, make an acceptable contribution for our conjecture.
The estimate is clearly best possible in the case of linear varieties.

Our first substantial result shows that it is those points that lie 
on lines in the hypersurface $X$ which are the only real difficulty. 
In fact we shall be able to tackle a rather more general counting 
problem in which
we consider points in a box, rather than a cube.  To describe this
situation we let $B_i \geq 1$ for $1\leq i \leq n$, and write
$\ma{B}=(B_1,\ldots,B_n)$ and 
\beq\lab{V}
V=\prod_{i=1}^n B_i.
\eeq
We then define
\[
N(F;\ma{B})=\#\{\x \in Z^n: F(\x)=0, ~|x_i| \leq B_i, ~(1\leq i \leq
n)\}.
\]
In particular we clearly have $N(F;\ma{B})=N(F;B)$ in the notation 
above, whenever $B_i=B$ for $1 \leq i \leq n$. 
Moreover we shall define $N_1(F;\ma{B})$
to be the number of points counted by
$N(F;\ma{B})$, but which do not lie on any line contained in $X$.
In the special case $B_i=B$ for  $1\leq i \leq n$, we 
merely write $N_1(F;B)$ for  $N_1(F;\ma{B})$.  Finally we write
$\|F\|$ for the maximum modulus of the coefficients of
$F$.  We then have the following two results, which will be 
established in \S \ref{deduction}.

\begin{thm}\lab{nolines}
Let $F\in \Z[x_1,\ldots,x_n]$ be an irreducible form of degree $d$.  
Then we have
$$
N_1(F;\ma{B})\ll V^{(n-2)/n}(\log \|F\|V)^{(n-2)(2n^2+n+3)/3},
$$
and
$$
N_1(F;B)=O(B^{n-2+\ve}).
$$
\end{thm}

\begin{cor}\label{nolinesatall}
Let $F\in \Z[x_1,\ldots,x_n]$ be an irreducible form of degree $d$,
and suppose that the variety $X$ is not a union of lines.  Then
$$
N(F;B)=O(B^{n-2+\ve}).
$$
\end{cor}

Handling the lines on the variety $X$ requires considerable further
work, and we have only been successful in extending the results proved in
\cite{annal} for $n=3$ and $n=4$, to the case $n=5$ corresponding
to threefolds.  We shall establish the following result, which provides 
further evidence for the main conjecture.

\begin{thm}\lab{3fold}
Let $F \in \Z[x_1,\ldots,x_5]$ be an irreducible form of degree $d
\geq 2$.  Then we have
$$
N(F;B)=O(B^{3+\ve}).
$$
\end{thm}

In fact our attack on the conjecture shows that every
point counted by $N(F;\ma{B})$ must lie on one of a small number of 
linear subspaces of
the hypersurface $X$.  The contribution from the points lying on
linear spaces of dimension $0$ (that is to say, from individual
points) will be satisfactory from the point
of view of the conjecture, while those lying on linear spaces of
higher dimension will be problematic.
In order to state the result we must first introduce some
more notation.  
Let $\sfl \subseteq \Z^n$ be any integer lattice of dimension $m\ge 2$ and
determinant $\det \sfl$, and define the set
\[
S(F;B,\sfl)=\{\x \in \sfl\cap Z^n: F(\x)=0, 
~\mbox{$\max_{i} |x_i|  \leq B$}\}. 
\]
Let $s_1 \leq \cdots \leq s_m$
be the successive minima of $\sfl$ with respect to the Euclidean
length.   
A brief discussion of successive minima will be included in \S \ref{pre}.
We then have the following result.

\begin{thm}\lab{main1}
Let $F \in \Z[x_1,\ldots,x_n]$ be a form of degree
$d$, not necessarily irreducible, but not vanishing identically on
$\sfl$. Let $p_0$ be the smallest prime $p_0>n$, and  let $B\geq
s_m$. Then there exist lattices $\sfm_1, \ldots, \sfm_J \subseteq \sfl$, 
such that the linear space $M_j$ spanned by each $\sfm_j$ lies in the
variety $X$, and such that
\beq\lab{res}
J \ll 
\left(\frac{B^m}{\det \sfl}\right)^{(m-2)/m} (\log
\|F\|B)^{(m-2)(2m^2+m+3)/3}.
\eeq
Moreover the successive minima of $\sfm_j$ are all
$O(B(\log\|F\|B)^{m(m+1)})$ and
\[S(F;B,\sfl)\subseteq\bigcup_{j=1}^J\sfm_j.\]
Associated to each proper sublattice $\sfm_j \subset \sfl$ of 
dimension $m-h$, are positive integers $q_1,\ldots,q_h$ and a 
vector $\ma{w}\in\Z^n$, all depending on
$\sfm_j$, and such that $q_i$ is a power of a prime $p_i\equiv
m+1-i\mod{p_0}$.  These have the property that
\beq\lab{hb1}
\sfm_j\subseteq\{\x\in \sfl: ~\x\equiv \rho\w \mod{q_1\cdots q_h}~
\mbox{some $\rho\in\Z$}\},
\eeq
and
\beq\lab{hb2}
\det\sfm_j\ll(\det\sfl)\prod_{i=1}^h\frac{q_i^{m-(i-1)}}{B}.
\eeq
\end{thm}

Evidently we will have $\dim \sfm_j=1+\dim M_j$, since we are using
the dimension of $M_j$ as a projective space.  The above result allows
us to provide information about the heights of the linear
spaces $M_j$ as follows.  This is crucial for the proof of Theorem
\ref{3fold}.  The height of a linear space will be defined in \S
\ref{cor-heights}.

\begin{thm}\lab{heights-planes}
For any $B \geq 1$ and $n \geq 4$, 
let $F \in \Z[x_1,\ldots,x_n]$ be an irreducible non-zero form of degree
$d$.   Let $\sfl=\Z^n$ in Theorem \ref{main1} and suppose that  
$$
M_1, \ldots, M_J \subseteq X 
$$ 
are the resulting linear spaces, which cover the points 
of $S(F;B,\Z^n)$.  Let $1\leq j \leq J.$  Then $M_j$ is defined over
$\Q$, and whenever the dimension
of $M_j$ is non-zero we have
$$
H(M_j) \ll B(\log||F||B)^{n^3}.
$$
\end{thm}

We stress that Theorem \ref{main1} refers to an 
arbitrary integer lattice $\sfl$, while Theorem \ref{heights-planes}
considers only $\Z^n$.  By change of basis, we shall see in \S 
\ref{deduction} that there is an approximate 
equivalence between points of a general lattice in a cube, and points 
of $\Z^n$ in an arbitrary box.  However in order to estimate 
$H(M_j)$ in Theorem \ref{heights-planes} it is necessary to restrict 
attention to points of $\Z^n$ in a cube.

It is interesting to see what Theorem \ref{nolines} has to say
about the quantity $N(F;\ma{B})$ in
the case of irreducible curves of degree $d \geq 2$.
Taking $n=3$ we have
$$
N(F;\ma{B})\ll V^{1/3} (\log \|F\|V)^{8}.
$$
This may be compared with work of the second author \cite[Theorem
3]{annal}, which implies that   
\beq\lab{d+1} 
N(F;\ma{B})\ll V^{1/(d+1)+\ve}. 
\eeq 
To see this we suppose for the moment that $B_1 \geq B_2 \geq B_3$, 
and note that an irreducible ternary form must 
contain a monomial in which $x_1$ appears explicitly and also a 
monomial in which $x_3$ does not appear.  Hence we may take 
$$
T \geq \max\{B_1B_3^{d-1},B_2^d\}\geq V^{d/(d+1)} 
$$ 
in the notation of \cite[Equation (1.8) of 
Theorem 3]{annal}, which suffices to establish (\ref{d+1}). 
Thus it follows that a direct application of Theorem \ref{nolines} is 
weaker than this 
earlier result, except possibly when $d=2$.  However, in this latter
case we are able to provide the following cleaner result.

\begin{thm}\lab{quadric}
Let $Q\in \Z[x_1,x_2,x_3]$ be a non-singular quadratic form.
Then we have
$$
N(Q;\ma{B})\ll V^{1/3}.
$$
\end{thm}

In
fact it is easy to see that this is  best possible whenever $B_i=B$ for
$1\leq i \leq 3$.
As a direct consequence of Theorem \ref{quadric} we are able to deduce
a further corollary, 
which sharpens \cite[Theorem  $2$]{cubic} and \cite[Corollary $2$]{annal}. 

\begin{cor}\lab{quadric'}
Let $Q \in \Z[x_1,x_2,x_3]$ be a non-singular quadratic form with
matrix $\ma{M}$.  Let $\D=|\det \ma{M}|$ and write $\D_0$ for the the 
highest common factor of the $2\times 2$ minors of $\ma{M}.$  Then 
$$
N(Q;\ma{B}) \ll \left\{ 1+ \left(\frac{V\D_0^{3/2}}{\D}\right)^{1/3} 
\right\}d(\D). 
$$
\end{cor}

This is superior to \cite[Theorem  $2$]{cubic} in three respects. 
Firstly we have replaced the exponent $1/3+\ve$ in   
\cite[Corollary $2$]{annal} by $1/3$.   
This is a direct result of our Theorem  
\ref{quadric}.  Secondly we have replaced $\D_0^2$ by $\D_0^{3/2}$. 
However it is already 
explicit in the proof of  \cite[Theorem  $2$]{cubic} that 
this is permissible.  Finally we have replaced 
the function $d_3(\D)$ by $d(\D)$. 
It is apparent from a closer inspection of the proof of 
\cite[Theorem  $2$]{cubic} that 
this too is allowable.

Theorem \ref{nolines} can also be used to good effect in the case $n=4$
of surfaces.  We assume without loss of generality that $B_1\geq B_2
\geq B_3 \geq B_4$, and suppose that 
$F \in \Z[x_1,\ldots,x_4]$ is any form of degree $d$.  Then Theorem
\ref{nolines} implies that
$$
N_1(F;\ma{B})\ll V^{1/2} (\log \|F\|V)^{26}.
$$
It is natural to ask, as in the case of curves above, whether or not we 
might expect a similar upper
bound to hold for the quantity $N(F;\ma{B})$ if $F$ is
irreducible  and has degree $d \geq 2$.  However simple examples of the type 
$$
F=x_3F_1+x_4F_2,
$$
for suitable forms $F_1,F_2 \in \Z[x_1,\ldots,x_4]$ of degree $d-1$, 
demonstrate that
we may have $N(F;\ma{B}) \gg B_1B_2$ whenever the surface $F=0$
contains the line $x_3=x_4=0$.  It is then a trivial matter to
deduce that we only have 
$V^{1/2} \geq B_1B_2$ when $B_1=B_2=B_3=B_4.$  These remarks show that
the following result, which will be established in \S \ref{pf-surface}, is
essentially best possible. 

\begin{thm}\lab{surface}
Let $F\in \Z[x_1,\ldots,x_4]$ be an irreducible form of degree $d \geq
2$ and suppose that $B_1\geq B_2 \geq B_3 \geq B_4$.
Then we have
$$
N(F;\ma{B})\ll (B_1B_2)^{1+\ve}.
$$
\end{thm}

{\bf Acknowledgement.} While working on this paper, the first
author was supported by
EPSRC grant number GR/R93155/01.

\section{Preliminaries}\lab{pre}

Our arguments will require a number of facts from the geometry of numbers,
which it will be convenient to collect together in this section. 
Most of these results can be found in the second author's work
\cite[\S 2]{square}, for example.

Let $\sfl\subseteq\Z^{n}\subset\R^{n}$ be a lattice of dimension
(or rank) $r$, and let 
$\ma{b}_{1},\ldots,\ma{b}_{r}$ be any 
basis for $\sfl$. Then we define the determinant $\det \sfl$ of
$\sfl$, to be the
$r$-dimensional volume of the parallelepiped generated by
$\ma{b}_{1},\ldots,\ma{b}_{r}$. This is independent of the choice of 
basis.  In fact if $\ma{M}$ denotes the $n\times r$ matrix formed from 
the vectors $\ma{b}_1,\ldots,\ma{b}_r$, then 
$\ma{M}^T\ma{M}$ is a real positive symmetric matrix and 
we have 
\beq\lab{basis-det} 
(\det \sfl)^2= \det(\ma{M}^T \ma{M}). 
\eeq 
We say that $\sfl$ is primitive if it is not properly
contained in any other $r$-dimensional sublattice of $\Z^{n}$.

In addition to these facts,  we shall also make repeated use of the
successive minima of a lattice.
Let $\sfl\subseteq\Z^{n}$ be a lattice of dimension $r$,
and let $\langle \ma{a}_1, \ldots, \ma{a}_r\rangle$ denote the
$\Z$-linear span of any set of vectors $\ma{a}_1, \ldots,
\ma{a}_r \in \R^n$.
Then we construct a minimal basis of $\sfl$ in the following manner.
First let $\ma{m}_1 \in \sfl$ be any non-zero vector for which the Euclidean
length $|\ma{m}_1|$ is least.  Next let $\ma{m}_2 \in
\sfl\setminus{\langle \ma{m}_1\rangle}$ be
any vector for which $|\ma{m}_2|$ is least.  Continuing in this way
we obtain a basis $\ma{m}_1, \ldots, \ma{m}_r$ for $\sfl$ in which 
$|\ma{m}_1| \leq \cdots \leq |\ma{m}_r|$.
The successive minima of $\sfl$ with respect to the Euclidean length
are merely the numbers $s_i=|\ma{m}_i|,$ for $1 \leq i \leq r$.
Then we have 
\beq\lab{rog6}
\prod_{i=1}^r s_i \ll \det \sfl \le \prod_{i=1}^r s_i. 
\eeq
In fact the upper bound in (\ref{rog6}) 
is a special case of the more general 
inequality 
$\det \sfl \le \prod_{i=1}^r |\ma{b}_i|,$ 
which holds for any basis $\ma{b}_{1},\ldots,\ma{b}_{r}$ for $\sfl$. 
A fundamental property of successive minima and minimal
bases is recorded in the following result, for which see Davenport 
\cite[Lemma 5]{Dav}.

\begin{lem}\lab{basis}
Let $\sfl\subseteq\Z^{n}$ be a lattice of dimension $r$, with
successive minima $s_1\leq \cdots \leq s_r$.
Then $\sfl$ has a basis $\ma{m}_1,\ldots,\ma{m}_r$ such
that $|\ma{m}_i|=s_i$ for $1 \leq i \leq r$, and with the property that
whenever one writes $\x\in\sfl$ as 
$$
\x=\sum_{i=1}^{r}\lambda_{i}\ma{m}_i,
$$ 
then
$\lambda_{i}\ll  |\x|/s_i,$
for $1 \leq i \leq r.$
\end{lem}

We shall also need a result which describes how the successive minima
of a lattice relate to the successive minima of any sublattice.

\begin{lem}\lab{mono}
Let $\sfl' \subseteq \sfl\subseteq\Z^{n}$ be lattices
of dimension $r$.
Suppose that $\sfl'$ has successive minima $s_1'\leq \cdots \leq
s_r'$ and that $\sfl$ has successive minima $s_1\leq \cdots \leq
s_r$.  Then we have 
$$
s_i' \geq s_i, \qquad (1 \leq i \leq r),
$$
and $\det \sfl' \gg \det \sfl.$
\end{lem}
\begin{proof}
To prove the first part of the lemma we suppose for a contradiction that
$$
s_1'\leq \cdots \leq s_{i}' < s_i,
$$ 
for some $1 \leq i \leq r$.
Now let $\ma{m}_1',\ldots, \ma{m}_r' \in \sfl'$ be the basis vectors
described in Lemma \ref{basis}, and analogously for
$\ma{m}_1,\ldots, \ma{m}_r \in \sfl$.  
Then it follows from 
the construction of the successive minima that the $i$ linearly
independent vectors $\ma{m}_1', \ldots, \ma{m}_i'$ must
belong to the $\Z$-linear span of $\ma{m}_1, \ldots,
\ma{m}_{i-1}$.  This is clearly impossible, which thereby
establishes that 
$s_i' \geq s_i$ for $1 \leq i \leq r.$

Finally we note that the second part of the lemma follows from the
first part and (\ref{rog6}).
\end{proof}

Next we recall a result of the second author \cite[Proof of Theorem
4]{annal}.  We shall say that a non-zero form defined over $\Z$ is 
primitive if the
highest common factor of its coefficients is $1$.

\begin{lem}\lab{theta}
Let $F \in \Z[x_1,\ldots,x_n]$ be a primitive form of degree
$d \geq 2$.  Then either
$$
\|F\| \ll B^{\theta}, \quad \theta= d\left(\colt{n-1 + d}{n-1}\right),
$$
or else  there exists a form 
$G \in \Z[x_1,\ldots,x_n]$ of degree $d$, not proportional to $F$,
such that every $\x \in  S(F;B,\Z^n)$ also satisfies the
equation $G(\x)=0$. 
\end{lem}

We shall also need a number of basic results about the counting
function for points on varieties.
Let $N \geq 3$ and let $Y \subset \bfP^{N-1}$ be  
a variety. Recall that for any $\x\in\Z^N$ we write $[\x]$ for the 
corresponding point in $\bfP^{N-1}(\Q)$.  We define the height of any 
rational point  $x \in \bfP^{N-1}(\Q)$ to be 
$$
H(x)=\max_{1\leq i \leq N}|x_i|,
$$
provided that $\x=(x_1,\ldots,x_{N}) \in Z^{N}$ satisfies $x=[\x]$. 
We shall write $Y(\Q)$ for $Y\cap \bfP^{N-1}(\Q)$, and  
$$
N_Y(B)=\# \{x \in Y(\Q): H(x)\leq B\}. 
$$
Thus in the case of Corollary \ref{nolinesatall}, for example, 
we have  
$N_X(B)=\frac{1}{2}N(F;B)$, since $\x$ and $-\x$ represent the same
point in projective space.

Our next result estimates $N_M(B)$, 
for  an $m$-dimensional linear space $M$ in $\bfP^{N-1}$.

\begin{lem}\lab{planecount}
Let $M \subseteq
\bfP^{N-1}$ be an $m$-dimensional linear space
 defined over $\Q$. Write 
\[
\sfm=\{\x\in\Z^N:[\x]\in M\}\cup\{\ma{0}\}, 
\]
so that $\sfm$ is a lattice of dimension $m+1$, and let
$s_1,\ldots,s_{m+1}$ be the successive minima of $\sfm$.  
Then
\[
N_M(B) \ll \prod_{j=1}^{m+1}\left(1+\frac{B}{s_j}\right). 
\]
Moreover if $m=1$ then
\[N_M(B)\ll 1+\frac{B^2}{\det\sfm}.\]
\end{lem}

\begin{proof}
Choose a basis $\ma{m}_1,\ldots,\ma{m}_{m+1}$ 
for $\sfm$ as in Lemma \ref{basis}.  Then 
if $\x\in \sfm$ with $\max_i|x_i|\le B$, we will have
\[
\x=\sum_{j=1}^{m+1}\la_j\ma{m}_j 
\]
with $\la_j \ll B/s_j$.  Thus the number of such points is 
\[
\ll \prod_{j=1}^{m+1}\left(1+\frac{B}{s_j}\right), 
\]
as required.

When $m=1$ we see that
\[
\x=\la_1\ma{m}_1+\la_2\ma{m}_2 
\]
with $\la_j\ll B/s_j$.  When $s_2\gg B$ this implies that $\la_2=0$, 
and since $N_M(B)$ counts only primitive vectors $\x$ we must have
$\x=\pm\ma{m}_1$. Thus $N_M(B)\ll 1$ when $s_2\gg B$.
On the other hand, when $s_1 \le s_2\ll B$ we find that 
\[ N_M(B)\ll
\left(1+\frac{B}{s_1}\right)\left(1+\frac{B}{s_2}\right)\ll 
1+\frac{B^2}{s_1s_2}\ll 1+\frac{B^2}{\det\sfm},\]
by (\ref{rog6}).  This suffices for the lemma.
\end{proof}

We shall also need the following estimates for points on curves and surfaces.
\begin{lem}\lab{curvecount}
Let $C$ be an irreducible curve in $\bfP^{N-1}$ of degree 
$D \geq 2$.  Then 
\[
N_C(B)\ll_{D,N}B^{2/D+\ve}\ll_{D,N}B^{1+\ve}. 
\]
Similarly, if $S \subset \bfP^{N-1}$ is any irreducible surface of degree
$D \geq 2$, we have 
\[N_S(B)\ll_{D,N}B^{2+\ve}.\] 
\end{lem}
The first assertion 
may be found in Broberg \cite[Corollary 1]{brob1}, for example, while
the second is in the first author's work \cite[Lemma 1]{browning}.

We end this section by establishing Theorem \ref{trivthm}, for which 
we let $\hat Y \subset \A^{N}$ be the  
affine cone above $Y$. Then $\hat{Y}$ is an 
irreducible affine variety
of degree $D$ and dimension
$m+1$.  Now for any irreducible affine variety $T \subset \A^{\nu}$ of 
degree $\del$ and
dimension $\mu$, we define the quantity
$$
M_T(B)=\#\{\ma{t} \in \Z^{\nu}: \ma{t} \in T, ~\mbox{$\max_{i} |t_i|
  \leq B$}\},
$$
and proceed by proving that 
\beq\lab{triv-affine}
M_T(B)=O_{\del,\nu}(B^\mu).
\eeq
This will clearly suffice to establish Theorem \ref{trivthm}, since 
$$
N_Y(B) \leq M_{\hat{Y}}(B). 
$$
We shall establish (\ref{triv-affine}) by induction on $\mu$.  
Since an irreducible variety of dimension zero contains just one 
point, the estimate is trivial whenever $\mu=0$.  Assume that $\mu
\geq 1$.  Since $T$ is irreducible we may find an index $1 \leq a\leq \nu$ such
that $T$ intersects the hyperplane
$t_a=\al$ properly, for any $\al \in \C$. Let $H_\al$ denote this
hyperplane.  We thereby obtain the upper bound
$$
M_T(B) \leq \sum_{|\al| \leq B} M_{T\cap H_\al}(B).
$$
Since $T\cap H_\al$ has dimension at most $\mu-1$ for every $\al$, 
and decomposes into at most $D$ irreducible exponents, 
an application of the  induction hypothesis implies that  
$M_{T\cap H_\al}(B)=O_{\del,\nu}(B^{\mu-1})$.  This suffices to complete the proof of
(\ref{triv-affine}), and so completes the proof of Theorem \ref{trivthm}.

\section{Deduction of Theorem \ref{nolines} and Corollary
\ref{nolinesatall}}\lab{deduction} 

In this section we use Theorem \ref{main1} to deduce Theorem
\ref{nolines} and Corollary  \ref{nolinesatall}.  
Recall the definition (\ref{V}) of $V$ and let
$$
b_i=\left\lfloor V/B_i\right\rfloor,
$$
where $\lfloor \alpha \rfloor$ denotes the integer part of any $\alpha \in \R$.
We define the lattice 
$$ 
\sfg(\ma{B};n) = \{\x \in \Z^{n}: b_i \mid x_i, ~(1\leq i \leq n)\}.
$$
In particular $\sfg(\ma{B};n)$ has rank $n$, successive minima
$b_1,\ldots,b_n$ (though not necessarily in that order) and determinant
\beq\lab{fall}
\det \sfg(\ma{B};n)  = \prod_{i=1}^{n} \left\lfloor V/B_i\right\rfloor 
\gg V^{n-1}.
\eeq
With these definitions in mind, we see that
under the linear transformation $\rom{Diag}(b_1,\ldots,b_{n})$, 
the image of the region $|x_i| \leq B_i$ lies inside 
a cube  with sides at most $2V$.  This therefore establishes the
useful inequality
$$ 
N(F;\ma{B}) \ll \#S(F;V,\sfg(\ma{B},n)),
$$
allowing us to move from considering points contained in a lopsided
region to points contained in a suitable cube.
In order to apply Theorem \ref{main1}, we observe that
the largest successive minimum $s_n$ of $\sfg(\ma{B};n)$ satisfies
$s_n \leq V$.   Indeed, if we assume that $B_1 \geq \cdots \geq B_n$
say, then $s_i=b_i$ for $1 \leq i \leq n$.
Now if $\x$ lies on some linear space $M \subseteq X$
of dimension at least $1$, then it trivially lies on some projective
line contained in $X$, and so is not to be counted by
$N_1(F;\ma{B})$. 
Taking $B=V$ and $\sfl=\sfg(\ma{B};n)$ in the statement of Theorem
\ref{main1}, we therefore deduce the first assertion of Theorem
\ref{nolines}, via (\ref{fall}).  

In order to deduce the second part of Theorem \ref{nolines} 
we employ Lemma \ref{theta}.
Since we may assume that $F$ is primitive, we deduce that
either $\log \|F\| \ll \log B$,  in which case we are done, or else 
that there exists a
non-proportional form $G$, of degree $d$, 
such that every $\x \in  S(F;B,\Z^n)$ lies
on the intersection $F=G=0$.   This is a union of $O(1)$ 
irreducible varieties of codimension 2 in $\bfP^{n-1}$, and degree at 
most $d^2$.  In this case the required estimate therefore follows from 
Theorem \ref{trivthm}.
This completes the proof of Theorem \ref{nolines}.

To deduce Corollary \ref{nolinesatall} we observe that it suffices,
after Theorem \ref{nolines}, to control the number of points lying on
lines contained in $X$.  
For any positive integer $m\leq n-1$ we let 
$\Grass(m,n-1)$ denote the variety parameterising $m$-dimensional
linear subspaces of $\bfP^{n-1}$.   
It is well-known that $\Grass(m,n-1)$ can be embedded in $\bfP^{N-1}$
via the Pl\"ucker embedding, where $N=\binom{n}{m+1}$.
Furthermore for any irreducible
hypersurface $Y$ in $\bfP^{n-1}$ we let 
$$
F_m(Y)=\{M \in \Grass(m,n-1): M \subset Y\}
$$
denote the corresponding Fano variety of $m$-dimensional
linear spaces contained in $Y$. 
Here we follow a common abuse of
notation and identify linear spaces $M \subset \bfP^{n-1}$ of
dimension $m$ with the corresponding point $M \in \Grass(m,n-1)$.
We now record the following basic result, which allows us to control
the degrees of several varieties relating to $F_m(Y)$, and which will
also be used in \S\S \ref{threefolds},\ref{pf-3fold} below.

\begin{lem}\lab{deg-fm}
Let $Y \subset \bfP^{n-1}$ be an irreducible hypersurface of degree
$D$ and let $y \in Y$. Then the three varieties 
$$
F_m(Y), \quad \bigcup_{M \in F_m(Y)}M, \quad \{M \in F_m(Y): y \in
M\},
$$ 
all have degrees which can be bounded in terms of $D$ and $n$ alone.
\end{lem}

\begin{proof}
That the degree of $F_m(Y)$ may be bounded in terms of $D$ and $n$ alone
follows by using the defining equation for $Y$ to write down the explicit 
equations for $F_m(Y)$.  To see the second part it will suffice to establish
that the degree of $\bigcup_{M \in \Psi}M$ does not exceed the
degree of $\Psi$, for any irreducible component $\Psi \subseteq F_m(Y)$. 
But this follows from a straightforward modification to \cite[Example
19.11]{harris}, much along the lines of \cite[Lemma 2.3.2]{broberg}.
Indeed this latter result establishes precisely this inequality in the
case $m=\dim \Psi=1$.

Finally we assume without loss of generality that 
$y=[1,0,\ldots,0]$, and let 
$p_{i_1\cdots i_{m+1}}$ be the Pl\"ucker coordinates of
$\Grass(m,n-1)$, 
for $1 \leq i_1<\cdots <i_{m+1} \leq n$.
Then it follows that the $m$-dimensional linear spaces in
$\bfP^{n-1}$ which pass through $y$  are parameterised by the
hyperplane
$$
\Lambda: \quad p_{i_1\cdots i_{m+1}}=0, \qquad (1 < i_1<\cdots <i_{m+1} \leq n),
$$
in $\Grass(m,n-1)$.  Hence 
$$
\{M \in F_m(Y): y \in
M\}=\Lambda \cap F_m(Y), 
$$
which clearly suffices to complete the proof of the lemma.
\end{proof}

We can now complete the proof of Corollary \ref{nolinesatall}.   
By hypothesis we may assume that $\bigcup_{L \in F_1(X)}L$ is a proper
subvariety of $X$.  Moreover Lemma \ref{deg-fm} ensures that it is  
a union of $O(1)$ irreducible varieties, each of degree 
$O(1)$. The result therefore follows from Theorem \ref{trivthm}.

\section{Conics}\lab{conics}

In this section we establish Theorem \ref{quadric}.
We may assume at the outset that $Q$ is primitive, since we can always
remove any factor common to the coefficients in the equation $Q=0$.
Let $\D$ be the discriminant of $Q$  and fix a choice of $r \in \N$.   
(In fact $r=109$ will suffice.)  By Bertrand's postulate we can choose primes
$p_1,\ldots,p_r$, with
\beq\lab{bertrand}
cV^{1/3} \leq p_1 <\cdots <p_r \ll_r V^{1/3},
\eeq
where $c$ is an absolute constant to be chosen in due course.
Now either there exists $1 \leq i\leq r$ for which $p_i \nmid \D$, or else
$$
|\D| \geq \prod_{i=1}^r p_i \gg V^{r/3}.
$$

We begin by disposing of the latter case.  Since $|\Delta|\le 6||Q||^3$ it
will follow that 
$$
||Q||\gg V^{r/9}\gg (\max B_i)^{r/9}.  
$$
We now apply 
Lemma \ref{theta} with $d=2$ and $n=3$, so that $\theta=12$.  This
shows that if $r>108$ then
there exists a ternary quadratic form $R$, not proportional to $Q$,
such that every $\x$ counted by $N(Q;\ma{B})$ also satisfies the
equation $R(\x)=0$.   But then B\'ezout's theorem reveals that  
$ N(Q;\ma{B}) \leq 4,$  
which is satisfactory.

We may now concentrate on the case in which there is a prime $p$ in
the range $cV^{1/3}\le p\ll V^{1/3}$ given by (\ref{bertrand}), with
the property that $p\nmid\D$.  Our argument now depends on the
following lemma.  We state the result in more generality than is needed for the
proof of Theorem \ref{quadric}, so that it may be applied later in the
treatment of Theorem \ref{main1}.

\begin{lem}\lab{M_k}
Let 
$H(x_1,\ldots,x_m)\in\Z[x_1,\ldots,x_m]$ be a form of degree $d$.  Let
$q$ be a power of a prime, and let 
$\x\in\Z^m$ be a vector for which $q|H(\x)$ and such that
\beq\lab{rog1}
{\rm h.c.f}\Big(q\,,\,\frac{\partial H(\x)}{\partial x_1}\,,\,\ldots\,,
\,\frac{\partial H(\x)}{\partial x_m}\Big)=1. 
\eeq
Then there is at least one vector $\x^{(1)}\in\Z^m$ satisfying both
(\ref{rog1}) and $q^2|H(\x^{(1)})$, and for which $\x^{(1)}\equiv\x\mod{q}$.
Write
$$
\sfm=\sfm(\x,q,H)=
\{\w\in \Z^m: ~\w\equiv \rho\x \mod{q}~\mbox{some $\rho\in\Z$}, ~q^2
\mid \w.\nabla H(\x^{(1)})\}.
$$
Then $\sfm$ is independent of the choice of $\x^{(1)}$.  Moreover
$\sfm$ is a lattice of dimension
$m$ and determinant $\det \sfm = q^m.$
Finally, if $\ma{t}\in Z^m$ satisfies $H(\ma{t})=0$, and if there exists
$\lambda\in\Z$ for which $\ma{t}\equiv\lambda\x\mod{q}$, then
$\ma{t}\in\sfm$. 
\end{lem}

\begin{proof}
The existence of a suitable $\x^{(1)}$ is an immediate consequence of
Hensel's Lemma.  Specifically, on writing $\x^{(1)}=\x+q\y^{(1)}$, we find that
$q^2|H(\x^{(1)})$ if and only if $\y^{(1)}.\nabla H(\x)\equiv
-q^{-1}H(\x)\mod{q}$, and this is always solvable for $\y^{(1)}$, 
by (\ref{rog1}).
Now suppose that $\x^{(2)}=\x+q\y^{(2)}$, with $\y^{(2)}.\nabla H(\x)\equiv
-q^{-1}H(\x)\mod{q}$, and moreover that $\w\equiv \rho\x \mod{q}$.  
To show that $\sfm$ is
independent of the choice of $\x^{(1)}$ it will suffice to demonstrate
that $q^2\mid \w.\nabla H(\x^{(1)})$ if and only if
$q^2\mid \w.\nabla H(\x^{(2)})$.  However
\begin{align*} 
\w.\nabla H(\x^{(i)})&=\w.\nabla H(\x+q\y^{(i)})\\
& \equiv\w.\nabla H(\x)+q\sum_{j,k=1}^nw_jy_k^{(i)}
\frac{\partial^2 H}{\partial x_j \partial x_k}(\x)\mod{q^2}\\
& \equiv\w.\nabla H(\x)+q\sum_{j,k=1}^n\rho x_jy_k^{(i)}
\frac{\partial^2 H}{\partial x_j \partial x_k}(\x)\mod{q^2}\\
& \equiv\w.\nabla H(\x)+q(d-1)\rho \y^{(i)}.\nabla H(\x)\mod{q^2}\\
&\equiv\w.\nabla H(\x)-(d-1)\rho H(\x)\mod{q^2}
\end{align*}
for $i=1,2$.  It follows that $q^2\mid \w.\nabla H(\x^{(1)})$
if and only if $q^2\mid \w.\nabla H(\x^{(2)})$, as required.

It is trivial that $\sfm$ is a lattice, since it is clearly closed 
under addition.  Moreover if $\y\in\Z^m$ then $q^2\y\in\sfm$ (taking 
$\rho=0$), whence $\sfm$ must have dimension $m$.  To compute
$\det\sfm$ we observe that if we put 
$\w=\rho\x^{(1)}+q\ma{z}$, then we have 
\begin{align*} 
\w.\nabla H(\x^{(1)})&= \rho\x^{(1)}.\nabla H(\x^{(1)})+ 
q\ma{z}.\nabla H(\x^{(1)})\\ 
&= \rho dH(\x^{(1)})+q\ma{z}.\nabla H(\x^{(1)}) \\ 
&\equiv  q\ma{z}.\nabla H(\x^{(1)})\mod{q^2},
\end{align*}
whence 
$$
\sfm= 
\{\w\in \Z^m: ~\w=\rho\x^{(1)}+q\ma{z}~\mbox{some
  $\rho\in\Z,\;\ma{z}\in\Z^m$}, ~q\mid \ma{z}.\nabla H(\x^{(1)})\}.
$$
We now observe that we have an inclusion of lattices $q^2\Z^m \subseteq 
\sfm \subseteq \Z^m$. Hence in order to calculate the determinant $\det 
\sfm$, which is equal to the index of $\sfm$ in $\Z^m$ as an additive subgroup,
it will suffice to calculate the index $[\sfm: q^2\Z^m]$. 
Indeed we then have
\beq\lab{index}
\det \sfm = [\Z^m: \sfm] = \frac{[\Z^m: q^2\Z^m]}{[\sfm: q^2\Z^m]} 
=\frac{q^{2m}}{[\sfm: q^2\Z^m]}. 
\eeq
We begin by considering the cosets of $\sfm$ modulo $q^2\Z^m$.
In view of the coprimality constraint (\ref{rog1}) there are $q^{m-1}$
possible values for $\ma{z}$ modulo $q$ satisfying $q\mid \ma{z}.\nabla
H(\x^{(1)})$. Moreover there are $q^2$ possible values for $\rho$ modulo
$q^2$.  Finally, each value of $\w$ may be decomposed as 
$\w\equiv\rho\x^{(1)}+q\ma{z}\mod{q^2}$ in $q$ ways.  It follows that
$\sfm(\x,q,H)$ has exactly $q^m$ cosets modulo $q^2\Z^m$, and so (\ref{index})
implies that $\det\sfm=q^{2m}q^{-m}=q^m$, as required.

For the final part of the lemma we note that if $\ma{t}\equiv
\lambda\x\mod{q}$, then $\ma{t}\equiv\lambda\x^{(1)}\mod{q}$.  Thus
there is a vector $\ma{z}\in\Z^m$ such that
$\ma{t}=\lambda\x^{(1)}+q\ma{z}$. It follows that
\[ 0=H(\ma{t})\equiv \lambda^d H(\x^{(1)})+
q\lambda^{d-1}\ma{z}.\nabla H(\x^{(1)})\equiv
q\lambda^{d-1}\ma{z}.\nabla H(\x^{(1)})\mod{q^2}.\]
Moreover we must have $\hcf(\lambda,q)=1$ since $\ma{t}$ 
is primitive and $\ma{t}\equiv\lambda\x\mod{q}$.  This allows us to
conclude that $q|\ma{z}.\nabla H(\x^{(1)})$.  We must
then have
\begin{align*}
\ma{t}.\nabla H(\x^{(1)})&=\lambda\x^{(1)}.\nabla H(\x^{(1)})+
q\ma{z}.\nabla H(\x^{(1)})\\
&=\lambda dH(\x^{(1)})+
q\ma{z}.\nabla H(\x^{(1)})\\
&\equiv 0\mod{q^2},
\end{align*}
so that $\ma{t}\in\sfm$.  This completes the proof of the lemma.
\end{proof}

We can now complete the proof of Theorem \ref{quadric}.  We have shown
that we may assume there is a prime $p\nmid\D$ in the range
$cV^{1/3}\le p\ll V^{1/3}$.  The projective variety $Q(\x)=0$ has
exactly $p$ points over $\mathbb{F}_p$, and we aim to show that there
are at most $2$ points counted by $N(Q;\ma{B})$ lying above each one.
This will suffice for the theorem.  We therefore fix a vector $\x\in
Z^3$ with $Q(\x)\equiv 0 \mod{p}$, and note that (\ref{rog1}) holds 
for $\x$ with $q=p$ since $p \nmid \D$.  Next we 
count vectors $\w\in Z^3$ satisfying $Q(\w)=0$, 
\beq\lab{rog2}
|w_1|\le B_1, \quad |w_2|\le B_2, \quad |w_3|\le B_3. 
\eeq
and such that there exists $\rho\in\Z$ with $\w\equiv \rho\x
\mod{p}$.  According to Lemma \ref{M_k} we will have
$\w\in\sfm=\sfm(\x,p,Q)$. We now consider the map
\[\phi(y_1,y_2,y_3)=(B_2B_3y_1\,,\,B_1B_3y_2\,,\,B_1B_2y_3).\]
Then $\phi(\sfm)$ has determinant $V^2\det\sfm=V^2p^3$.  Let $s_1\le
s_2\le s_3$ be the successive minima of $\phi(\sfm)$, and let
$\ma{m}_1,\ma{m}_2,\ma{m}_3$ be the corresponding basis 
vectors described in Lemma \ref{basis}. 
It then follows from (\ref{rog6}) that 
$s_3\geq V^{2/3}p\gg V$.  On the other hand, if 
(\ref{rog2}) holds, then $\phi(\w)$ has Euclidean length $\ll V$.  
Hence Lemma \ref{basis} implies that, if we choose the absolute  
constant $c$ in (\ref{bertrand}) to be sufficiently large,
the vector $\phi(\w)$ must lie in the 2-dimensional lattice spanned by
$\ma{m}_1$ and $\ma{m}_2$.  Thus, given $\x$, all 
corresponding points $\w$ must lie not only on the conic $Q(\w)=0$ but
also on a certain line.  There are therefore at most 2 such 
points $\w\in Z^3$.  This completes the proof of Theorem \ref{quadric}.

\section{Proof of Theorem \ref{main1} --- variable moduli}\label{VM}

The purpose of this section is to examine how the ideas of the
previous section may be adapted to hypersurfaces of dimension 2 or more.
It is clear that \S \ref{conics} depends crucially on condition
(\ref{rog1}).  For any prime $p$, any lattice $\sfl \subseteq \Z^n$ of 
dimension $m \geq 2$, and any form $G\in\Z[x_1,\ldots,x_n]$ of degree $d$, 
we therefore introduce the set
\[
S(G;B,\sfl,p)=\{\x \in S(G;B,\sfl): p\nmid\nabla G(\x)\}.
\]
We could then mimic the argument from \S \ref{conics} to show that if
$p\gg B$ then the
points in $S(G;B,\sfl,p)$ lie in $O(B^{m-2})$ sublattices $\sfl'$ of $\sfl$,
each of dimension $m-1$.  In \S \ref{conics}, where $G=Q$ is a 
non-singular ternary quadratic form and $\sfl=\Z^3$, each sublattice 
$\sfl'$ corresponds to a line, and each line contains at most $2$ 
points on the conic $Q=0$.  
When $m= 4$, for example, the situation is more complicated.
Each sublattice $\sfl'$ corresponds to a plane in $\bfP^3$, so that the
corresponding points of $G=0$ typically lie on a curve.  In order to
estimate the number of such points we would want to know the size of
$\det \sfl'$, and in order to get a good estimate we would want
$\det\sfl'$ to be large.  Since it is clear that we produce different
lattices $\sfl'$ for different points $\x\mod{p}$ we must anticipate
difficulties in trying to make $\det\sfl'$ large for all $\x\mod{p}$.
We shall overcome these problems
by using moduli which are powers $p^k$ of $p$,
and allowing different values of $k$ for different vectors $\x$.

It is now time to state the result we shall use.

\begin{lem}\label{VMlemma}
Let $\sfl\subseteq\Z^n$ be a lattice of dimension $m\ge 2$,  
with largest successive minimum at most $B$. 
Let $p$ be a prime with $p\nmid\det\sfl$.
Then there is an integer $I\ll p^{m(m-2)}(B^m/\det\sfl)^{(m-2)/m}$ and lattices
$\sfl_1,\ldots,\sfl_I\subset\sfl$ of dimension $m-1$, such that
\[S(G;B,\sfl,p)\subseteq\bigcup_{i=1}^I\sfl_i.\]
For any $i$, the successive minima of $\sfl_i$ are all $O(p^mB)$.
Moreover for any $\al\geq 1$ we have 
\[ 
\#\{i:\det\sfl_i\le \alpha\det\sfl\}\ll p^{m-2}(\alpha B)^{(m-2)/m}.
\]
\end{lem}

\begin{proof}
Let $\ma{e}_1,\ldots,\ma{e}_m$ be a basis for $\sfl$ and set
\[H(y_1,\ldots,y_m)=G(y_1\ma{e}_1+\cdots+y_m\ma{e}_m).\]
Let $\phi:\sfl\rightarrow\Z^m$ be the map
\[\phi:\x=y_1\ma{e}_1+\cdots+y_m\ma{e}_m\mapsto \y.\]
Then
\[\frac{\partial H}{\partial y_j}=\ma{e}_j.\nabla G(\x),\]
and since $p\nmid\det\sfl$ we see that 
$p\nmid\nabla G(\x)$ implies $p\nmid\nabla H(\y)$.  

We proceed by defining a finite tree for each non-singular projective
point $\y$ on the mod $p$ reduction of $H=0$.  There will be
$O(p^{m-2})$ such points, uniformly in $H$. We define the 
``depth'' of the root node to be 1, and if a node $N'$ is an immediate
successor to a node $N$ we define the depth of $N'$ to be 1 more than
the depth of $N$.  This terminology may differ slightly from that
employed by other authors.

We shall define an
equivalence relation on the set
\[
\{\ma{v}\in(\Z/p^k\Z)^m: {\rm h.c.f.}(v_1,\ldots,v_m,p)=1\}, 
\]
by taking $\ma{v}\sim_k\ma{v}'$ if and only if $\ma{v}=\lambda\ma{v}'$
for some unit $\lambda$ of $\Z/p^k\Z$.  The tree we construct will
then have each node of depth $k$ labelled by an equivalence class
$[\ma{v}]_k$, such that, if $\ma{v}=\ma{w}+p^k\Z^m$ then
$p^k|H(\ma{w})$ and $p\nmid\nabla H(\ma{w})$.  The root node
is thus labelled by the equivalence class $[\y]_1$, where $\y$ is the 
point mentioned above. 

To a node with label $[\ma{v}]_k$, where $\ma{v}=\ma{w}+p^k\Z^m$, 
we associate the lattice $\sfn=\phi^{-1}\sfm(\ma{w},p^k,H)$, 
in the notation of
Lemma \ref{M_k}.  Clearly the lattice will be independent of the
choice of $\ma{w}$.  
We then declare the node to be a leaf if the largest
successive minimum $s_m$ of this lattice satisfies $s_m>cB$.  Here 
$c$ is an absolute constant whose value will become clear in due 
course. 

We must now compute the determinant of $\sfn$.  If we take
$\y_1,\ldots,\y_m\in \Z^m$ to be column vectors forming a basis of 
$\sfm(\ma{w},p^k,H)$, then $\sfn$ will have a basis
$\ma{E}\y_1,\ldots,\ma{E}\y_m$, where $\ma{E}$ is the $n\times m$
matrix with columns $\ma{e}_1,\ldots,\ma{e}_m$.  Then
(\ref{basis-det}) implies that 
\[
(\det \sfn)^2=\det(\y_i^T\ma{E}^T\ma{E}\y_j) 
\]
and
\[
(\det\sfl)^2=\det(\ma{E}^T\ma{E}). 
\]
Since $\ma{E}^T\ma{E}$ is a real positive symmetric matrix, we may 
write $\ma{E}^T\ma{E}=\ma{F}^T\ma{F}$ for some real $m\times m$ matrix 
$F$.  Thus if $\ma{Y}$ is the $m\times m$ matrix with columns
$\y_1,\ldots,\y_m$ we will have
\begin{align*}
(\det \sfn)^2&=\det(\y_i^T\ma{E}^T\ma{E}\y_j)\\ 
&=\det(\y_i^T\ma{F}^T\ma{F}\y_j)\\
&=\det(\ma{Y}^T\ma{F}^T\ma{F}\ma{Y})\\
&=(\det \ma{Y})^2(\det \ma{F})^2, 
\end{align*}
and
$$ 
(\det \sfl)^2 =\det(\ma{E}^T\ma{E}) 
=\det(\ma{F}^T\ma{F}) =(\det \ma{F})^2.  
$$
Finally, 
\[
\det \sfm(\ma{w},p^k,H)=|\det(\ma{Y})|, 
\]
whence
\beq\lab{5.0}
\det\sfn=(\det \sfm(\ma{w},p^k,H))(\det \sfl)=p^{km}\det\sfl, 
\eeq
by Lemma \ref{M_k}.

Thus the largest successive minimum of $\sfn$ satisfies
\beq\label{rog4}
s_m\geq (\det\sfn)^{1/m}= p^k(\det\sfl)^{1/m}, 
\eeq
by (\ref{rog6}).
For each node $[\ma{v}]_k$ which is not a leaf, we take its immediate
successors to be those nodes labelled by $[\ma{u}]_{k+1}$, for which
$\ma{u}=\ma{w}+p^{k+1}\Z^{m}\in(\Z/p^{k+1}\Z)^{m}$ satisfies both
$\ma{v}=\ma{w}+p^k\Z^m$ and $p^{k+1}|H(\ma{w})$.  In particular we
have $p \nmid \nabla H(\ma{w})$ since $[\ma{v}]_k$ is a node.
We claim that there will be precisely $p^{m-2}$ such immediate successors.
This is a simple application of Hensel's Lemma.  Clearly it will 
suffice to show that for any node $[\ma{v}]_k$, there are $p^{m-1}$ 
vectors $\ma{u}=\ma{w}+p^{k+1}\Z^{m}$ for which 
$\ma{v}=\ma{w}+p^k\Z^m$ and $p^{k+1}|H(\ma{w})$.
Suppose that $\ma{v}=\ma{w}+p^k\Z^m$, with $p \nmid \nabla H(\ma{w})$
and $H(\ma{w})=p^{k}w$, say.
Then we count values of $\ma{u} \in (\Z/p^{k+1}\Z)^{m}$ for which
$\ma{u}\equiv \ma{w} \mod{p^{k}}$ and $p^{k+1} \mid H(\ma{u})$.
Writing $\ma{u}=\ma{w}+p^{k}\ma{z}$, we deduce that
$$
H(\ma{u}) \equiv H(\ma{w}) + p^{k} \ma{z}.\nabla H (\ma{w}) \mod{p^{k+1}},
$$
and so seek values of $\ma{z}$ for which 
$$
w+\ma{z}.\nabla H (\ma{w}) \equiv 0 \mod{p}.
$$
But this has exactly $p^{m-1}$ solutions $\ma{z} \mod{p}$, since 
$p \nmid \nabla H(\ma{w})$.
Therefore there are $p^{m-2}$ immediate successors
above any given node which is not a leaf.

We will need to establish an upper bound for $s_m$.  In fact we shall
show that
\beq\lab{s_m}
s_m \ll p^m B.
\eeq
Suppose first that the root node $[\ma{y}]_1$ is a leaf, and let
$\sfn=\phi^{-1}\sfm(\ma{y},p,H)$ be the corresponding lattice.  Then upon 
observing that
we have an inclusion of lattices 
$p^2\sfl\subseteq\sfn$, it follows from Lemma \ref{mono}
that the successive minima of $\sfn$ are at most $p^2$
times as large as those of $\sfl$.  Since the largest successive
minimum of $\sfl$ is at most $B$ by hypothesis, we therefore  
have $s_m\le p^2B$, which is satisfactory for (\ref{s_m}).
Whenever the root node is not a leaf we may consider 
the node $[\ma{v}']_{k-1}$, say, 
immediately preceding the leaf $[\ma{v}]_k$ for some $k \geq 2$.  
Let the  lattices corresponding to $[\ma{v}]_k$ and $[\ma{v}']_{k-1}$
be $\sfn$ and $\sfn'$ respectively, with successive minima
$s_1,\ldots,s_m$ and $s_1',\ldots,s_m'$.  
We note that $\sfn\subseteq\sfn'$, as one sees directly
from the definition in Lemma \ref{M_k}.
But then it follows from Lemma \ref{mono} that
$$
s_i \geq s_i', \qquad (1 \leq i \leq m).
$$
Since $[\ma{v}']_{k-1}$ is not a leaf we have $s'_{m}\ll B$.  
Therefore (\ref{rog6}) yields 
\begin{align*}
s_m&\ll \frac{\det \sfn}{s_1\cdots s_{m-1}}\\
&\leq \frac{\det \sfn}{s'_1\cdots s'_{m-1}}\\
&\leq \frac{\det \sfn}{\det \sfn' /s'_m }, 
\end{align*}
whence an application of (\ref{5.0}) completes the proof of (\ref{s_m}).

We may conclude from (\ref{rog4}) and (\ref{s_m}) that 
$$
p^{k-m} \ll \left(\frac{B^m}{\det \sfl}\right)^{1/m}. 
$$
In view of this we now see that 
the tree is finite, with total depth
$$
k_0\ll 1+\frac{\log(B^m/\det\sfl)}{\log p},
$$
and has at most
$$
p^{(m-2)(k_0-1)}\ll
p^{(m-1)(m-2)}\left(\frac{B^m}{\det\sfl}\right)^{(m-2)/m} 
$$ 
leaf nodes.  The reader should also recall that we have $O(p^{m-2})$ such
trees to consider, giving $O(p^{m(m-2)}(B^m/\det\sfl)^{(m-2)/m})$ 
leaf nodes in total.

Now suppose that $\x\in S(G;B,\sfl,p)$.  Then $\y=\phi(\x)$ satisfies
$H(\y)=0$ and $p\nmid \nabla H(\y)$.  
The construction of the various trees is such that
there is then a leaf node $[\ma{v}]_k$ in one of the trees with the
property that
$\y\equiv\lambda\ma{v}\mod{p^k}$, for some $\lambda$.  It follows from
Lemma \ref{M_k} that $\y$ belongs to the corresponding lattice 
$\sfm(\ma{v},p^k,H)$, 
and hence that $\x$ belongs to the corresponding lattice $\sfn$.
As we have seen in Lemma \ref{basis}, there exists a basis
$\ma{m}_{1},\ldots, \ma{m}_{m}$ of the lattice $\sfn$ such 
that if one writes any $\x\in \sfn$ as $\x=\sum_{i=1}^m
\la_i\ma{m}_{i}$, then  $\la_i \ll |\x|/s_i$.  Thus if  
$\max_i|x_i|\le B$ and 
$s_m>cB$ with $c$ a suitably chosen constant, we must have $\lambda_m=0.$
Hence if $\x\in S(G;B,\sfl,p)$ corresponds to the leaf node
$[\ma{v}]_k$ then $\x$ must belong to the lattice $\sfl_i$, say, 
spanned by $\ma{m}_{1},\ldots, \ma{m}_{m-1}$.  The number of such 
lattices $\sfl_i$ is the total number of leaf nodes, which is 
$O(p^{m(m-2)}(B^m/\det\sfl)^{(m-2)/m})$.

We next consider the successive minima of
$\sfl_i$.  Clearly these are no larger than the corresponding values of
$s_m$, which by (\ref{s_m}) are $O(p^mB)$.  Indeed it is clear that
any set of linearly independent vectors in $\sfl_i$ is also a set of
linearly independent vectors in $\sfn$, whence the successive minima
of $\sfl_i$ must be precisely $s_1,\ldots,s_{m-1}$.  According to 
(\ref{rog6}) and (\ref{5.0}) we therefore have
\[
\det\sfl_i\gg\prod_{j=1}^{m-1}s_j=s_m^{-1}\prod_{j=1}^m s_j\geq 
\frac{\det\sfn}{s_m}\gg 
p^{km-m}B^{-1}\det\sfl. 
\]
Thus $\det\sfl_i\le\alpha\det\sfl$ implies that $p^{k}\ll p(\alpha
B)^{1/m}$.  Since each tree has at most $p^{(m-2)(k-1)}$ leaf nodes of
depth $k$ it follows that there are at most $O(p^{m-2}(\alpha
B)^{(m-2)/m})$ possible values for $i$.  This completes the proof of
the lemma.
\end{proof}

For future use we also note that
\beq\lab{future}
\det\sfl_i\leq\prod_{j=1}^{m-1}s_j=s_m^{-1}\prod_{j=1}^m s_j\ll 
\frac{\det\sfn}{s_m}\ll p^{km}B^{-1}\det\sfl, 
\eeq
by a further application of (\ref{rog6}), (\ref{5.0}) and the fact 
that $s_m \gg B$ by construction.

\section{Proof of Theorem \ref{main1} --- the induction}\label{IA}

In this section we shall complete the proof of Theorem \ref{main1}
using an induction argument.  For this we shall fix $n$, and use
induction on $m$.  We take $m=2$ as the base for the
induction.  In this case $S(F;B,\sfl)$ consists of $O(1)$ points $\ma{t}$ 
with $|\ma{t}|\ll B$.  Thus 
we can take $J\ll 1$, and each $\sfm_j$ 
will be spanned by the corresponding $\ma{t}$.  Moreover, each
$\sfm_j$ will have dimension $1$ and successive minimum $O(B)$, so 
that (\ref{hb1}) holds with $\ma{w}=\ma{t}$ for any $q_1$.   
Choosing a sufficiently large prime $q_1\equiv 2\mod{p_0}$ will
therefore suffice to ensure that (\ref{hb2}) holds.
Hence the required
results hold for the base case $m=2$ and 
we may assume henceforth  that  $m \geq 3$.

Our first task is to show how Lemma \ref{VMlemma} may be applied, and
this is achieved via the following result.

\begin{lem}\lab{splitting}
Let $\mathcal{P} = \log^2(\|F\| B)$.  Then
there is an integer $t\ll 1$ and forms 
$G_1,\ldots,G_t\in\Z[x_1,\ldots,x_n]$, with degrees at
most $d$, not vanishing identically on $\sfl$.  Moreover
\[
||G_j||\ll ||F||, \qquad (1 \leq j \leq t),  
\]
and
\beq\lab{cup1}
S(F;B,\sfl)\subseteq\bigcup_{1 \leq j \leq t} \;  
\bigcup_{\mathcal{P}\leq p\ll\mathcal{P}}S(G_j;B,\sfl,p). 
\eeq
Here the primes $p$ are restricted to satisfy the conditions $p\equiv
m\mod{p_0}$ and $p\nmid\det\sfl$. 
\end{lem}

\begin{proof}
We begin by defining
$$
\mcal{G}=\left\{\left.\frac{\partial^{a_1+\cdots+a_n}F}
{\partial^{a_1}x_1\cdots\partial^{a_n}x_n}\right|_{\sfl}\not\equiv 0:
\;  (a_1,\ldots,a_n) \in \Z_{\geq 0}^n\right\}.
$$ 
Thus $\mcal{G}$ is just the set of all partial derivatives of $F$ 
which do not vanish identically on $\sfl$. 
We then take the
forms $G_j$ to be the forms in the set $\mathcal{G}$.  Thus it 
suffices to show (\ref{cup1}), the remaining claims being obvious.  

Let 
$\x  \in S(F;B,\sfl)$  and consider the subset
$$
\mcal{G}(\x)=\{G \in \mcal{G}: G(\x)=0\}.
$$
In particular we observe that $\mcal{G}(\x)$ is non-empty since 
$F \in \mcal{G}$.
We take $G=G_{\x}$ to be any form in the set $\mcal{G}(\x)$ which has
minimal degree, and proceed to show that  $\x$ is a non-singular 
point of $G$.  To see this
we first note that $G(\x)=0$, since $G \in \mcal{G}(\x)$.  Moreover, 
our choice of
$\mcal{G}$ implies that $G$ does not vanish identically on $\sfl$.
Hence $G$ must have degree at least $1$, since $G(\x)=0$.  
But then it follows that the
components of $\nabla G$ cannot all vanish identically on $\sfl$, 
and so there
exists $1 \leq i\leq n$ for which the partial derivative $\partial
G/\partial x_i$ does not vanish identically on $\sfl$.  Finally we deduce from
the minimality of $\deg G$ that 
$$
\frac{\partial G}{\partial x_i}(\x) \neq 0,
$$
from which it follows that $\nabla G(\x) \neq \ma{0}$ as claimed.  
It is perhaps instructive to remark at this point that if $F$ is a 
non-singular form, as it was in \S \ref{conics}, then we can take
$t=1$ and $G_1=F$. In fact our construction of 
auxiliary forms $G_j$ is used purely in order to handle the 
contribution to  
$S(F;B,\sfl)$ from the singular locus of $F=0$.

Now let $G=G_\x$ be the form constructed above, with 
$\nabla G(\x)\not=\ma{0}$, and let $p$ be the least prime number 
$p \geq \mathcal{P}$ for which $p \equiv m \mod{p_0}$ and 
$$
p \nmid (\det\sfl)\nabla G(\x).
$$
We now observe that 
$$
\det\sfl\leq \prod_{j=1}^ms_j\leq B^m, 
$$
by (\ref{rog6}) and the hypotheses of Theorem \ref{main1}, whence
$\log\det\sfl\ll\log B$.
Recalling that $G$ is non-singular at $\x$,
a trivial modification to the  proof of
\cite[Lemma $4$]{annal} then shows that $p$ exists,
and that
$$
\mathcal{P}\leq p \ll \mathcal{P}.
$$
It follows that $\x$ belongs to the corresponding set $S(G;B,\sfl,p)$,
and this suffices for the lemma.
\end{proof}

We are now ready to begin the proof of Theorem \ref{main1}, using
induction on $m$. According to Lemmas \ref{VMlemma} and 
\ref{splitting} we 
have
\[
S(F;B,\sfl)\subseteq\bigcup_{1 \leq j \leq t} 
\bigcup_{\mathcal{P}\leq p\ll\mathcal{P}}\bigcup_{i=1}^I\sfl_i 
\]
for suitable lattices $\sfl_1, \ldots, \sfl_I \subset\sfl$ of 
dimension $m-1$. 
Moreover the successive minima of these lattices are all 
\beq\lab{rog7}
\le c\mathcal{P}^mB=B_0,  
\eeq
say, for a suitable constant $c$ depending only on $d$ and $n$.  Hence 
\[
S(F;B,\sfl)\subseteq\bigcup_{\sfm} S(F;B_0,\sfm), 
\]
where $\sfm$ runs over all the lattices $\sfl_i$ as the forms $G_j$
and the primes $p$ vary.
The number of lattices $\sfm$ with $\det\sfm\le\alpha\det\sfl$ is
$$ 
O(\mathcal{P}^{m-1}(\alpha B)^{(m-2)/m}), 
$$ 
since there are $O(\mcal{P})$ primes overall. Similarly, one finds 
that the total number of 
lattices $\sfm$ is $O(\mcp^{m(m-2)+1}(B^m/\det\sfl)^{(m-2)/m})$. 
It follows from (\ref{rog7}) that
\beq\lab{rog9}
\det\sfm\ll \mcp^{m(m-1)}B^{m-1}.
\eeq

We now focus on a particular lattice $\sfm$. If $\sfm=\sfl_i$ arises from
$\sfm(\ma{w},p^k,H)$ in the proof of Lemma \ref{VMlemma}, then we take 
$q_1=p^k$.  In particular it follows that $q_1$ is a power of a prime 
$p_1=p\equiv m\mod{p_0}$.  
Our construction gives
\[\sfl_i\subset\sfn=\phi^{-1}\sfm(\ma{w},p^k,H),\]
from which it is clear that if $\x\in\sfl_i$ then
\beq\lab{q1}
\x\equiv\rho\phi^{-1}(\ma{w})\mod{q_1}
\eeq
for some $\rho\in\Z$.
If the form $F$ vanishes
on the lattice $\sfm$ then the variety $F=0$ contains the
corresponding linear space.  The number of such lattices is
satisfactory for (\ref{res}), and (\ref{q1}) suffices for (\ref{hb1}).
Moreover, (\ref{future}) suffices for (\ref{hb2}), since $h=1$.

We suppose from now on that $F$ does not vanish identically 
on $\sfm$.  We may therefore apply our induction hypothesis to
$S(F;B_0,\sfm)$ and conclude that it is covered by
\begin{align*}
&\ll \left(\frac{B_0^{m-1}}{\det \sfm}\right)^{(m-3)/(m-1)} 
(\log\|F\|B)^{(m-3)\{2(m-1)^2+(m-1)+3\}/3}\\
&\ll \left(\frac{B^{m-1}}{\det \sfm}\right)^{(m-3)/(m-1)}\mcp^{m(m-3)}
(\log\|F\|B)^{(m-3)\{2(m-1)^2+(m-1)+3\}/3}
\end{align*}
linear spaces, each contained in the variety $F=0$.  
It follows that the total number of linear
subvarieties arising from all
those lattices $\sfm$ for which 
\[\frac{1}{2}\alpha\det\sfl<\det\sfm\le\alpha\det\sfl\]
is
\begin{align*}
&\ll \mcp^{m-1}(\alpha B)^{(m-2)/m}\\
& \qquad \qquad\times 
\left(\frac{B^{m-1}}{\alpha\det \sfl}\right)^{(m-3)/(m-1)}\mcp^{m(m-3)} 
(\log\|F\|B)^{(2m^3-9m^2+13m-12)/3}\\
&\ll \frac{\alpha^{2/(m(m-1))}B^{(m^2-2m-2)/m}}{(\det\sfl)^{(m-3)/(m-1)}}
(\log||F||B)^{(2m^3-3m^2+m-18)/3}.
\end{align*}
In view of (\ref{rog9}) we must sum this over dyadic ranges for 
\[\alpha\ll \frac{B^{m-1}}{\det\sfl}\mcp^{m(m-1)},\]
producing a total 
\begin{align*}
&\ll\left(\frac{B^m}{\det\sfl}\right)^{(m-2)/m}
(\log||F||B)^{(2m^3-3m^2+m-6)/3},
\end{align*}
as required.  This establishes the bound (\ref{res}).

By our induction hypothesis, the successive minima of the lattices
that arise will be
\[
\ll B_0(\log ||F||B_0)^{(m-1)m}.  
\]
In view of our choice of $B_0$ this
is $O(B(\log||F||B)^{m(m+1)})$, as required.

According to our induction hypothesis $S(F;B_0,\sfm)$ is covered by
lattices $\sfm_j\subseteq\sfm$.  When $\sfm_j$ has dimension $(m-1)-h$
for $h \geq 1$,  there are associated to it positive integers
$q_1',\ldots,q_{h}'$,  and a vector $\ma{w}'\in\Z^n$ such that
\[
\sfm_j\subseteq\{\x\in \sfm: ~\x\equiv \rho\w' \mod{q_1'\cdots q_{h}'}~
\mbox{some $\rho\in\Z$}\}.
\]
Each $q_{i}'$ is a power of a prime $p_{i}'\equiv (m-1)+1-i\mod{p_0}$.
Now we have already noted in (\ref{q1}) that there exists an integer
$q_1$ which is a power of a prime $p_1 \equiv m \mod{p_0}$, such that
\[
\x\equiv\rho\phi^{-1}(\ma{w})\mod{q_1}
\]
for some $\rho \in \Z$, whenever $\x\in\sfm$.  For $1 \leq i \leq h+1$  
we define the integers
$$
q_i=\left\{
\begin{array}{ll}
q_1, & i=1,\\
q_{i-1}', & i>1.
\end{array}
\right.
$$
In particular it follows that each $q_i$ is a power of a prime $p_i
\equiv m+1-i \mod{p_0}$.  These congruence constraints ensure that $q_1$ is coprime to
$q_2\cdots q_{h+1}$.  Hence the Chinese Remainder Theorem
implies that there is a vector $\ma{w}\in\Z^n$ such that
\[\sfm_j\subseteq\{\x\in \sfm: ~\x\equiv \rho\w \mod{q_1q_2\cdots q_{h+1}}~
\mbox{some $\rho\in\Z$}\},\]
as required for (\ref{hb1}).

It remains to consider the bound (\ref{hb2}).  However on applying
(\ref{hb2}) for $\sfm$ we find that
\[
\det\sfm_j\ll(\det\sfm)\prod_{i=1}^{h}\frac{q_i'^{(m-1)-(i-1)}}{B_0} 
=(\det\sfm)\prod_{i=1}^{h}\frac{q_i'^{m-i}}{B_0}.
\]
Since $B_0\gg B$, the bound (\ref{future}) yields
\[
\det\sfm_j\ll q_1^mB^{-1}(\det\sfl)\prod_{i=2}^{h+1}
\frac{q_i^{m-(i-1)}}{B} 
\ll (\det \sfl)  
\prod_{i=1}^{h+1}\frac{q_i^{m-(i-1)}}{B},
\]
as required.  This completes the proof of Theorem \ref{main1}.

\section{Proof of Theorem \ref{heights-planes}}\lab{cor-heights}

We begin by defining the height of a linear space.  
For positive integers $m\leq n-1$, we recall from \S \ref{deduction}
the variety $\Grass(m,n-1)$ which parameterises $m$-dimensional linear
subspaces of $\bfP^{n-1}$, and which may be embedded in projective
space via the Pl\"ucker embedding.
Whenever $M \in \Grass(m,n-1)$ is defined over $\Q$, 
we define the height $H(M)$ of
$M$ to be the height of its coordinates in $\Grass(m,n-1)$, under the
Pl\"ucker embedding.  
Then for any  $M \in \Grass(m,n-1)$ which is defined over $\Q$, we let 
$$
\sfm=\{\x \in \Z^n: [\x] \in M\}\cup\{\ma{0}\}
$$
be the lattice associated to $M$.
It follows from \S \ref{pre} that there exists a basis
$\ma{e}_1,\ldots,\ma{e}_{m+1}$ of $\sfm$ such that 
$|\ma{e}_1|\cdots|\ma{e}_{m+1}|$ has the same order of magnitude as
$\det \sfm$.  In  fact it can be shown (see Schmidt \cite[Chapter I,
Corollary 5I]{schmidt}, for example) that 
\beq\lab{sch}
\det \sfm =H(M).
\eeq

The proof of Theorem \ref{heights-planes} is now
straightforward. Suppose that $\sfl=\Z^n$ in Theorem \ref{main1}, so
that $m=n$, and let $\sfm_j$ have dimension $l=n-h$.  Let the
successive minima of $\sfm_j$ be $s_1\le\cdots\le s_l$ and choose a basis
$\ma{e}_1,\ldots,\ma{e}_{l}$ for $\sfm_j$ as in Lemma \ref{basis}.  We now
observe, according to
(\ref{hb1}), that
there exists $\w\in\Z^n$ and
$\rho_1, \ldots, \rho_{l} \in\Z$ such that 
$$
\ma{e}_i \equiv \rho_i\ma{w} \mod{Q}, \qquad (1\leq i \leq l),
$$
where $Q=q_1\cdots q_h$.  Define $\sigma_i=1+\lfloor s_i\rfloor$ for 
$1\le i\le l$
and consider the lattice
$$
\sfn=\{(n_1\sigma_1,\ldots, n_{l}\sigma_{l}) \in \Z^{l}: 
~\rho_1n_1+\cdots+ \rho_{l}n_{l}\equiv 0 \mod{Q}\}.
$$
If ${\rm h.c.f.}(\rho_1,\ldots,\rho_l,Q)\neq 1$ then $\sfm_j$ contains 
no primitive vectors, and hence can make no contribution to 
$S(F;B,\Z^n)$.  We may therefore suppose that  
\[{\rm h.c.f.}(\rho_1,\ldots,\rho_l,Q)=1,\] 
whence $\sfn$ has rank $l$ and determinant 
\[\det\sfn=Q\prod_{i=1}^{l} \sigma_i \ll Q \det\sfm_j, \]
by (\ref{rog6}).  We write $t_1\le\cdots\le t_l$ for the successive
minima of $\sfn$ and choose a basis $\ma{c}^{(1)},\ldots,
\ma{c}^{(l)}$ for $\sfn$ as in Lemma \ref{basis}.  Thus in particular 
$|\ma{c}^{(j)}|=t_j$ for each $j$ and hence
\beq\lab{bb}
\prod_{j=1}^l|\ma{c}^{(j)}|=\prod_{j=1}^l t_j\ll\det\sfn\ll
Q\det\sfm_j,
\eeq
by a further application of (\ref{rog6}).  We now put
\[\ma{c}^{(j)}=(d_1^{(j)}\sigma_1,\ldots,d_l^{(j)}\sigma_l),\;\;\;(1\le
j\le l)\]
and
$$
\ma{f}_j= Q^{-1}(d_1^{(j)}\ma{e}_1+\cdots+d_l^{(j)}\ma{e}_l), \qquad
(1\leq j \leq l).
$$
By construction we have $\ma{f}_j\in\Z^n$ and 
$[\ma{f}_j] \in M_j(\Q)$. Moreover, since the vectors  
$\ma{c}^{(1)},\ldots,\ma{c}^{(l)}\in\Z^l$ 
are linearly independent it follows that the matrix
\[
\left(d_{i}^{(j)}\right)_{1\le i,j\le l}
\]
is invertible.  Bearing in mind that the vectors 
$\ma{e}^{(1)},\ldots,\ma{e}^{(l)}$ are linearly independent, it
follows that $\ma{f}^{(1)},\ldots,\ma{f}^{(l)}$ are also
linearly independent. 
Let
\[
\sfn^{(*)}=\langle \ma{f}^{(1)},\ldots,\ma{f}^{(l)}\rangle
\]
be the $\Z$-linear span of $ \ma{f}^{(1)},\ldots,\ma{f}^{(l)}$, 
and observe that $\det\sfn^{(*)} \ll
\prod_{j=1}^l|\ma{f}^{(j)}|$.  Now it is clear from the triangle
inequality that we have
$$
|\ma{f}^{(j)}|\leq Q^{-1}\sum_{i=1}^l|d_i^{(j)}||\ma{e}_i|,
$$
for each $1 \leq j \leq l$.   Employing (\ref{bb}) we therefore obtain
\begin{align*}
\det\sfn^{(*)} &\ll Q^{-l}\prod_{j=1}^l\left(\sum_{i=1}^l|d_i^{(j)}|s_i\right)\\
&\ll Q^{-l}\prod_{j=1}^l\left(\sum_{i=1}^l|d_i^{(j)}|\sigma_i\right)\\
&\ll Q^{-l}\prod_{j=1}^l|\ma{c}_i^{(j)}|\\
&\ll Q^{1-l}\det\sfm_j.
\end{align*}
Now write 
\[
\sfm_j^{(*)}=\{\x\in\Z^n:\,[\x]\in M_j\}\cup\{\ma{0}\} 
\]
for the lattice associated to the linear space $M_j$.  
It is worthwhile  highlighting  that although $\sfm_j$ spans $M_j$, it is
not necessarily the case  that $\sfm_j$ is equal to $\sfm_j^{(*)}$.
However we clearly have
\[
\sfn^{(*)}\subseteq\sfm_j^{(*)},
\]
so that
\[H(M_j)=\det \sfm_j^{(*)} \ll \det\sfn^{(*)}\ll Q^{1-l}\det\sfm_j,\] 
by (\ref{sch}) and Lemma \ref{mono}.
We now observe that (\ref{hb2}) yields
\[\det\sfm_j\ll\prod_{i=1}^h\frac{q_i^{n-(i-1)}}{B}\ll \frac{Q^n}{B^h}, \]
since we are taking $\sfl=\Z^n$.  Thus, on using the relation $l=n-h$
we obtain
\[H(M_j)\ll Q^{1+h}B^{-h}.\]
Alternatively, since all the successive minima of $\sfm_j$ are
$O(B(\log||F||B)^{n(n+1)})$, it follows from (\ref{rog6}) that
\[\det\sfm_j\ll B^{n-h}(\log||F||B)^{n(n+1)(n-h)},\]
so that
\[H(M_j)\ll Q^{1-l}B^{n-h}(\log||F||B)^{n(n+1)(n-h)}
\ll Q^{1-n+h}B^{n-h}(\log||F||B)^{n^3}.\]
We use the first of these two bounds for $Q\le B(\log||F||B)^{n^2}$,
and the second otherwise.  Theorem \ref{heights-planes} then follows.

\section{Proof of Theorem \ref{surface}}\lab{pf-surface}

Let $S \subset \bfP^3$ denote the surface $F=0$.  We begin by
considering the contribution to $N(F;\ma{B})$ arising from a line $L$,
defined over $\Q$ and lying in the surface $S$.  Let $\sfl$ be the integer
lattice corresponding to $L$. Then $H(L)=\det\sfl$, by
(\ref{sch}), and we claim that 
\begin{align}
\#\{\x\in Z^4 &: [\x] \in   L(\Q), ~|x_i| \leq B_i, ~(1\leq i \leq 4)\} 
\nonumber\\  
&=\#\{\x \in \sfl\cap Z^4: |x_i| \leq B_i, ~(1\leq i \leq 4)\} 
\nonumber\\  
&\ll 1+ \frac{B_1B_2}{\det \sfl}.\lab{tim1} 
\end{align} 
The box we are concerned with lies inside the ellipsoid
\[
\Big(\frac{x_1}{B_1}\Big)^2+\Big(\frac{x_2}{B_2}\Big)^2+\Big(\frac{x_3}{B_2}\Big)^2 
+\Big(\frac{x_4}{B_2}\Big)^2\le 4. 
\]
Arguing in $\mathbb{R}^4$, the intersection of this ellipsoid with any
2-dimensional linear space will be an ellipse of area $O(B_1B_2)$.
Thus, to establish (\ref{tim1}) it suffices to count points in
$\mathbb{R}^2$ lying in an ellipse, and which belong to a lattice
$\sfl'$ of determinant $\det{\sfl}$.  By using an appropriate
unimodular transformation $\ma{M}$ we may map this ellipse to a disc
of the same area.  We are therefore led to
count the number of lattice points contained in a disc of area $O(B_1B_2)$ which
lie on $\ma{M}\sfl'$.  Since $\det(\ma{M}\sfl')=\det(\sfl')=\det(\sfl)$,
an application of \cite[Lemma 1 (vi)]{annal} therefore yields the
bound (\ref{tim1}).

We now apply Lemma \ref{theta}.  If  
$\log||F||\gg\log B_1$ this shows that we 
may confine attention to points on a union of $O(1)$ curves $C$ in
$S$, each of degree $O(1)$.  For 
those curves $C$ of degree at least $2$ we use the bound
\[
N_C(B_1)\ll B_1^{1+\ve},
\]
which follows from Lemma \ref{curvecount}, while for the case in which
$C$ is a line we get a satisfactory estimate from (\ref{tim1}).
Henceforth we may therefore assume that $\log||F||\ll\log B_1$.

By the argument used to derive Theorem \ref{nolines} from
Theorem \ref{main1}, the points in
which we are interested lie on $O(V^{1/2}(\log B_1)^{26})$ linear subspaces
of $S$.  Those subspaces of dimension 0 are clearly satisfactory for Theorem
\ref{surface}, since $V^{1/2} \leq B_1B_2$.  Similarly those lines
$L\subset S$
which contain at most one rational point in the box under consideration
also make a satisfactory contribution.  The remaining lines have at
least two rational points in the relevant box, and hence have height
$O(B_1^2)$. 

Recall the definition of the Fano variety $F_1(S)=\{L \in \Grass(1,3):
L \subset S\}$ of lines in $S$.  Then Lemma \ref{deg-fm} implies that
the degree of $F_1(S)$ is bounded uniformly in terms of $d$.
Moreover, it is well-known and easy to prove that whenever $d \geq 2$,
we have 
$$
\dim F_1(S) \leq 1,
$$
and $F_1(S)$ contains no linear component of dimension 1.  In
particular it follows from Lemma \ref{curvecount} that
$$
N_{\Phi}(2H) \ll H^{1+\ve},
$$
for any irreducible component $\Phi \subseteq F_1(S)$.

Let $L\in F_1(S)$ be a line defined over $\Q$, with height $H< H(L) \leq 2H$ for
some $H\ll B_1^{2}$.  Then (\ref{tim1}) shows that each such line
contributes $O(1+B_1B_2/H)$ to $N(F;\ma{B})$.  We have just seen that
there are $O(H^{1+\ve})$ available lines.  However we also know that
the number of linear spaces $M_j$ is $O(V^{1/2}B_1^{\ve})
=O(B_1^{1+\ve}B_2)$.  It
follows that the total contribution from the  
lines $L=M_j \subset S$ which are defined over
$\Q$ and have height $H < H(L)\leq 2H$ is
$$
\ll B_1^{1+\ve}B_2+\frac{B_1B_2}{H}.H^{1+\ve} = 
B_1B_2(B_1^{\ve}+H^{\ve}). 
$$
Since we have $H(L)\ll  B_1^{2}$ we may sum this bound over dyadic
intervals with $H \ll B_1^{2}$ to obtain the overall estimate $\ll
(B_1B_2)^{1+2\ve}$.  
We complete the proof of Theorem \ref{surface} on
re-defining $\ve$.

\section{Geometry of cubic threefolds}\lab{threefolds}

In view of the results of Heath-Brown
\cite[Theorem 2]{annal} and of Broberg and Salberger \cite{broberg}
already cited it suffices for the proof of Theorem \ref{3fold}
to consider the situation in which $d=3$.
Thus we assume that the equation $F=0$ defines an irreducible cubic 
threefold $X \subset \bfP^4$, throughout this section. 
It is clear from Theorem \ref{main1} that we must examine the
possible lines and planes on cubic threefolds. 
For $m=1,2$, recall from \S \ref{deduction} the Fano variety 
$F_m(X)\subset \Grass(m,4)$   
of $m$-dimensional linear spaces contained in
$X$.  By Lemma \ref{deg-fm} it follows that $F_m(X)$ is an intersection 
of hypersurfaces whose
degree and number are bounded absolutely.

We first consider the dimension of $F_1(X)$, since the larger this
dimension is, the more difficult it will be to handle the contribution
from the lines on $X$.
Define the incidence correspondence
$$
I= \{(x,L) \in X \times F_1(X): x \in L\}.
$$
The fibre of $I$ over a point in $F_1(X)$ has dimension $1$. Hence it 
follows that the dimension of $I$ is $1+\dim F_1(X)$.
Now consider the projection onto $X$.  If 
$x \in X$ is a generic point then the tangent hyperplane
$\mathbb{T}_x(X)$ to $X$ at $x$ has dimension 3.  Thus the 
dimension of $X\cap \mathbb{T}_x(X)$ is $2$.  Moreover,  any
line $L \subset X$ which passes through $x$, must  be contained in
$\mathbb{T}_x(X)$. Hence the fibre of $I$ over a generic
point $x \in X$ has dimension at most $1$.  It follows that $\dim I 
\leq 4$, and so 
$$
 \dim F_1(X) \leq 3.
$$
In a precisely similar manner a consideration of
\beq\lab{J}
J= \{(x,P) \in X \times F_2(X): x \in P\}
\eeq
shows that
$$
 \dim F_2(X) \leq 1.
$$
We also record the following facts about $F_1(X)$ and $F_2(X)$.  

\begin{lem}\lab{degPhi}
Let $X \subset \bfP^4$ be an irreducible cubic threefold.  Then the 
following hold:
\begin{enumerate}
\item[(i)]
Suppose $\Psi\subseteq F_1(X)$ where $\Psi$ is a plane.  Then the lines
$L\in\Psi$ sweep out a plane in $X$.
\item[(ii)] $F_2(X)$ does not contain a line.
\end{enumerate}
\end{lem}
\begin{proof}
For the proof of part (i)
we observe that any plane in $\Grass(1,4)$ corresponds 
either to a
locus of lines passing through the same point, which are all contained
in the same hyperplane, or else to a locus of coplanar lines.  In the
first case we would find that $X$ contains a hyperplane, which is
impossible since $X$ is an irreducible cubic threefold.  In the
second case the lines sweep out a plane, as claimed.

Similarly for part (ii) we note that a line in $\Grass(2,4)$ corresponds 
to the set of planes of $\bfP^4$ containing a given line, and
contained in a fixed hyperplane.  They therefore sweep out the
hyperplane.  Since $X$ does not contain a hyperplane $F_2(X)$ cannot
contain a line. 
\end{proof}

Whenever the dimension of $F_1(X)$ is maximal we shall
call upon the following geometric result of Segre \cite{segre}.

\begin{lem}\lab{fano}
Let $X\subset \bfP^4$ be an irreducible cubic threefold, and suppose
that  $F_1(X)$ has dimension $3$. Then 
$$
\dim F_2(X)=1.
$$ 
\end{lem}

Recall that any plane contains a two dimensional family of lines.
Thus whenever $\dim F_2(X)=1$, we see that $X$ contains a three dimensional 
family of lines, all of which lie in
planes contained in  $X$.  In fact 
most of the lines in $X$ arise in this way, as is shown by 
the following result.

\begin{lem}\lab{l-i-p}
Let $X\subset \bfP^4$ be an irreducible cubic threefold for which
$F_2(X)$ has dimension $1$. 
Then every line in $X$ lies in a plane  
contained in $X$, with the possible exception of those lying on a
certain surface of degree $O(1)$ contained in $\Grass(1,4)$. 
\end{lem}

Suppose firstly that $X$ is a cone, and let $x \in  X$ be a vertex
point.  The family of lines in $X$ passing through $x$ is an algebraic
variety 
$$
\{L \in F_1(X): x\in L\}=Y_x, 
$$ 
say.  Moreover $Y_x$ has
dimension $2$ and degree $O(1)$ by Lemma \ref{deg-fm}. 
Let $L \in F_1(X)\setminus Y_x$.  
Then if we write $\overline{x,L}$ for the plane spanned by $x$ and $L$, we
clearly have
$$
L \subset \overline{x,L} \subset X.
$$
This establishes Lemma \ref{l-i-p} whenever $X$ is a cone.

We suppose henceforth that $X$ is not a cone. 
Throughout the proof of Lemma \ref{l-i-p}, we will use 
$H \in {\bfP^4}^*$ to denote a generic hyperplane.  Since
$X$ is not a cone, the hyperplane section 
$$
S_H=H \cap X
$$ 
is an irreducible
cubic surface which is not a cone.
We claim that $S_H$ contains infinitely
many lines.  Let $\Phi \subseteq F_2(X)$ be an irreducible component
of dimension $1$, and observe that 
\beq\lab{*1}
X = \bigcup_{P \in \Phi}P.
\eeq
In particular it clearly follows that 
$$
S_H= \bigcup_{P \in \Phi} (H\cap P).
$$
The intersection $H \cap P$ cannot ever be a plane, since $S_H$ is an
irreducible cubic surface.  Hence $H \cap P$ is always a line.  Thus
$S_H$ is a union of lines and hence contains infinitely many lines as claimed.

We now call upon the 
following classical result, which follows from the work of
Bruce and Wall \cite{b-w} for example.

\begin{lem}\lab{ruled}
Let $S \subset \bfP^3$ be a ruled  irreducible cubic surface.   Then
either $S$ is a cone or the singular locus of $S$ is a line.
\end{lem}

Let $Y \subset X$ denote the singular locus of $X$, and let $T_H
\subset S_H$ denote the singular locus of $S_H$.
Since we have already seen that $S_H$ contains infinitely many lines, 
Lemma \ref{ruled} implies that $T_H$ is a line.
An application of Bertini's theorem (in the form given by Harris 
\cite[Theorem $17.16$]{harris}, for example) therefore shows that
$$
H\cap Y =T_H \cong \bfP^1
$$
for $H \in {\bfP^4}^*$. 
It follows that $Y$ must be a plane.
Suppose that
$X$ is given by a cubic form $F(\x)=F(x_1,x_2,x_3,x_4,x_5)$.  After a
linear change of variables we may assume that the plane $Y$ is given
by $x_1=x_2=0.$  It follows that there exist quadratic forms $Q_1,Q_2$
such that
$$
F(\x)=x_1Q_1(\x)+x_2Q_2(\x).
$$
Upon considering the partial derivatives of $F$ with respect to
$x_1,x_2$ one deduces further that for $i=1,2$ the forms
$Q_i(0,0,x_3,x_4,x_5)$ must vanish identically, since $Y$ is a double plane.
Hence there exist linear forms $L_1,L_2,L_3$ such that
$$
F(\x)=x_1^2L_1(\x)+x_2^2L_2(\x)+x_1x_2L_3(\x).
$$
Moreover we may assume that $L_1,L_2,L_3$ are linearly independent,
since otherwise $X$ would be a cone.  We have therefore established
the following result, which may be of independent interest. 

\begin{lem}\lab{scroll}
Let $X \subset \bfP^4$ be an irreducible cubic threefold with $\dim F_2(X)=1$.
Then either $X$ is a cone or else $X$ takes the shape
$$
x_1x_2x_3+x_1^2x_4+x_2^2x_5=0.
$$
\end{lem}

We observe that for $X$ as in Lemma \ref{scroll}, there is a family of
planes given by 
\beq\lab{e-p}
\la x_1-\mu x_2= \la \mu x_3 +\mu^2 x_4+\la^2 x_5=0, 
\eeq
for any $[\la,\mu] \in \bfP^1$, in addition to the  plane
\beq\lab{e-p'}
x_1=x_2=0.
\eeq

It remains to consider the lines contained in $X$.  We hope to prove
that the generic such line is contained in one of the planes
(\ref{e-p}).  
Now any line  in $\bfP^4$ can be given parametrically by 
$$
\x=[a_1s+b_1t,a_2s+b_2t,a_3s+b_3t, a_4s+b_4t, a_5s+b_5t], 
$$
for suitable $a_i,b_i \in \overline{\Q}$.   
If $a_1b_2-a_2b_1 = 0$ the line is contained in a hyperplane $\la
x_1-\mu x_2=0$, and it is readily deduced that the line lies in one of
the planes (\ref{e-p}) or (\ref{e-p'}).  It therefore remains to
examine the case $a_1b_2-a_2b_1 \neq 0$, and here
it suffices to take $a_1=b_2=1$
and $a_2=b_1=0$.  By equating coefficients in the vanishing binary form 
$st(a_3s+b_3t)+s^2(a_4s+b_4t)+t^2(a_5s+b_5t)$, we conclude
that any line contained in $X$, which is not
contained in any of the planes (\ref{e-p}), must take the shape
$$
\al x_3=\be x_1+\gamma x_2, \quad \al x_4=-\be x_2, \quad \al
x_5=-\gamma x_1,
$$
for appropriate $[\al,\be,\gamma] \in \bfP^2$ such that $\al \neq 0$.
One readily finds that the family of all
such lines forms a surface in $\Grass(1,4)$.  Any line whose
Pl\"ucker cordinates do not lie on this surface will be contained in 
one of the planes
(\ref{e-p}) or (\ref{e-p'}).  This completes the proof of Lemma \ref{l-i-p}.

It is interesting to remark that
for $X$ as in Lemma \ref{scroll}, the Fano variety of planes $F_2(X)$ is
the union of a single twisted cubic  and an isolated point in
$\Grass(2,4)$. However we shall make no use of this fact in our work.
We end this section by considering points which lie on infinitely many
planes.

\begin{lem}\lab{infplanes}
Let $X$ be an irreducible cubic threefold.  If $X$ is a cone then its
set of vertices is either a single point or a line.  Moreover, if
$x\in X$ lies
on infinitely many planes in $X$, then $X$ is a cone with vertex 
$x$.    Indeed if $X$ is not a cone with vertex $x$, then  
$x$ lies on just $O(1)$ planes in $X$. 
\end{lem}

\begin{proof} 
If $x$ and $y$ are two distinct vertices of $X$ any point $z$ on 
the line $\overline{x,y}$
is also a vertex.  Thus the set of vertices, $V$ say, is necessarily a
linear space.  However if $V$ has dimension 2, and $p\in X\setminus
V$, then all of $\overline{p,V}$ is contained in $X$, which is 
impossible, since $X$ cannot contain a hyperplane.  
Now suppose that $x\in X$ lies
on infinitely many planes in $F_2(X)$.
Let $J$ be as defined in (\ref{J}), so that the fibre over
$x$ is infinite, and therefore 
of dimension 1 since $\dim F_2(X) \leq 1$.   
Let $\Phi$ be an irreducible component of this fibre, with dimension 1.  Then
(\ref{*1}) shows that $X$ is a cone with vertex $x$.  Thus if $x$ is
not a vertex of $X$ the fibre of $J$ over $x$, which may be written 
$$
\{P \in F_2(X): x \in P\}=Y_x,
$$
say, has dimension 0. In order to complete the proof of the lemma, it
therefore suffices to show that the degree of $Y_x$ is bounded
absolutely.  But this follows directly from Lemma \ref{deg-fm}.
\end{proof}

\section{Proof of Theorem \ref{3fold}}\lab{pf-3fold}

As already remarked we may restrict our attention to the case in which
$X$ has degree 3.  By combining Lemma \ref{theta} and Theorem \ref{trivthm} we
see that it suffices also to assume that $\log||F||\ll\log B$.  Then,
according to Theorem \ref{main1}, the points in
which we are interested lie on $O(B^3(\log B)^{58})$ linear subspaces
of $X$.  Subspaces of dimension 0 are clearly satisfactory for Theorem
\ref{3fold}, so it remains to consider the case in which $M_j$ is a
line or a plane.  In view of Theorem \ref{heights-planes} 
we know that $H(M_j)\ll B(\log B)^{125}$ in these cases.

Our argument begins by considering points which lie on certain planes
$P\in F_2(X)$.  We shall call such a plane ``$B$-good'' if it
contains three points $x_1,x_2,x_3$ in general position, with
$H(x_i)\le B(\log B)^{31}$.  We remark at once that whenever $M_j$ is
a plane it is automatically $B$-good if $B$ is large enough, since
Theorem \ref{main1} tells us that the successive minima of $\sfm_j$
are $O(B(\log B)^{30})$.  For any $B$-good plane $P$ we 
let $\sfm$ be the integer
lattice corresponding to $P$, so that $H(P)=\det\sfm$, by
(\ref{sch}).  Then the successive minima of $\sfm$ will be $O(B(\log
B)^{31})$, so that Lemma \ref{planecount} yields
$$
N_P(B)\ll \prod_{j=1}^{3}\left(1+\frac{B}{s_j}\right) 
\ll \prod_{j=1}^{3}\frac{B(\log B)^{31}}{s_j}, 
$$
since $B \leq B(\log B)^{31}$.  But then it follows that  
$$
N_P(B) \ll B^3(\log B)^{93}/\det\sfm = B^3(\log B)^{93}/H(P). 
$$

We proceed to consider the contribution from all $B$-good planes for
which $H<H(P)\le 2H$.  Consider those planes which belong to a
particular irreducible component $\Phi$ of $F_2(X)$. 
The number of such planes is at most $N_{\Phi}(2H)$.  Moreover $\Phi$ is
either a single point, or is a curve of degree at least 2, by Lemma
\ref{degPhi}. Since the degree of $\Phi$ is absolutely  bounded, by
Lemma \ref{deg-fm}, Lemma
\ref{curvecount} shows that 
$$ 
N_{\Phi}(2H)\ll H^{1+\ve}. 
$$  
It follows
that the total contribution from all $B$-good planes $P$ with height 
$H<H(P)\le 2H$ is
\[\ll \frac{B^3(\log B)^{93}}{H}.H^{1+\ve}.\]
Since we have $H(P)\le B^3(\log B)^{93}$ we may sum this bound over
dyadic intervals with $H\ll B^3(\log B)^{93}$ to obtain an estimate
$O(B^{3+2\ve})$ for the total contribution to $N(F;B)$ arising from
points on $B$-good planes.  On re-defining $\ve$ this is satisfactory
for Theorem \ref{3fold}.  It follows that we have a satisfactory
contribution whenever $M_j$ is a plane, and also whenever $M_j$ is a
line contained in a $B$-good plane.

We can now dispose of the case in which $X$ is a cone with a
line $L$ of vertices, such that $L$ contains two distinct rational points
$p,q$ with $H(p),H(q)\le B$.  In this instance, if $x\in X\setminus L$ is
a rational point with $H(x)\ll B$, then the plane $\overline{p,q,x}$ 
lies in $X$ and is $B$-good.  Thus the number of such points is
$O(B^{3+\ve})$ by the above.  Since the line $L$ contains  
$O(B^2)$ admissible
points, by Theorem \ref{trivthm}, we may conclude that Theorem
\ref{3fold} holds in this case.  Henceforth we shall assume that $X$
is not of the type just considered.  In the light of Lemma
\ref{infplanes} our assumption tells us that $X$ has at
most one vertex $v$ which is a rational point of height $H(v)\le B$. 

We now turn our attention to the case of lines $M_j$.  In view of our
treatment of $B$-good planes, it will be enough to consider lines that
are not contained in any $B$-good plane.  Let $L$ be a
line defined over $\Q$, and let $\sfl$ be the corresponding
2-dimensional integer lattice.  Fix a basis $\ma{m}_1, \ma{m}_2$ of 
$L$ of the type given 
in Lemma \ref{basis}.  Then we shall say that $L$ is ``$B$-exceptional'' if it
does not belong to a $B$-good plane, and if
either $\dim F_2(X)\not=1$, or if $\dim F_2(X)=1$ and $L$ corresponds
to one of the exceptional lines in Lemma  \ref{l-i-p},
or if $\dim F_2(X)=1$ but $X$ is a cone with vertex $[\ma{m}_1]$. 
We proceed by considering the contribution 
from $B$-exceptional lines $M_j$.  In each case we shall show that we have at
most a two dimensional family of lines.
If $\dim F_2(X)\not=1$ then $\dim F_1(X)\le 2$,
by Lemma \ref{fano}.  In the second case the claim is an
immediate consequence of Lemma \ref{l-i-p}.  In the third case
we note that the  set of lines passing through any given point of $X$ will have
dimension at most 2. The claim is therefore obvious when $X$ has at most one 
vertex.  If $X$ has a line of vertices, we have seen above that there is at  
most one such vertex $v$ which is defined over $\Q$ and which has height  
$H(v)\le B$.  However we 
observe that 
\[|\ma{m}_1|=s_1\le (s_1s_2)^{1/2}\ll  
(\det\sfl)^{1/2}=H(M_j)^{1/2}\ll B^{(1+\ve)/2}, \] 
by (\ref{rog6}), (\ref{sch}) and Theorem \ref{heights-planes}.  Thus 
$|\ma{m}_1|\le B$ if $B$ is large enough.  It follows that  
$[\ma{m}_1]=v$, so that $M_j$ 
must pass through the unique vertex $v$ of $X$ for which $H(v)\le B$. 
We therefore conclude in all cases that the set of $B$-exceptional 
lines $M_j$ has dimension at most 2.

Now it is clear that any $B$-exceptional line must be contained either
in $F_1(X)$ if $\dim F_2(X)\not=1$,
or to the closed set described in Lemma \ref{l-i-p} if 
$\dim F_2(X)=1$ and $L$ corresponds
to one of the exceptional lines described in this result, or to the set
$\{L \in F_1(X): [\ma{m}_1] \in L\}$ if $X$ is a cone with vertex $[\ma{m}_1]$. 
If $\Psi\subset F_1(X)$ is an irreducible component of any of these 
closed sets then we have
already seen that $\dim \Psi \leq 2$, and it follows from Lemma
\ref{deg-fm} and Lemma \ref{l-i-p} that $\Psi$ has degree $O(1)$.
Moreover, if $\Psi$ were  a plane then the $B$-exceptional lines 
$L\in \Psi$ would 
all lie in a plane $P$, say, contained in $X$, by Lemma \ref{degPhi}.  
The contribution from such components is therefore $O(B^3)$, by Theorem
\ref{trivthm},  and this is satisfactory for our theorem. Such components $\Psi$
may therefore be ignored.

We may now complete the treatment of the $B$-exceptional lines in much the same
way as we dealt with the $B$-good planes.  
If $L=M_j$ is a $B$-exceptional line then Theorem \ref{main1} implies that
the successive minima of the corresponding lattice $\sfm_j$ will be  
$O(B(\log B)^{30})$, so that Lemma \ref{planecount} yields
\begin{align*}
N_L(B) &\ll \prod_{j=1}^{2}\left( 1+\frac{B}{s_j}\right)\\
&\ll\prod_{j=1}^{2}\frac{B(\log B)^{30}}{s_j}\\
&\leq B^2(\log B)^{60}/\det\sfm_j\\ 
&= B^2(\log B)^{60}/H(L).
\end{align*}
We have now to count the number of $B$-exceptional lines $L$ for which
$H<H(L)\le 2H$.  We do this by considering those which belong to a
given subvariety $\Psi\subset F_1(X)$, where we may now restrict 
attention to the case in which either $\Psi$ has dimension at most 1,
or it has dimension
2 and degree at least 2.  For such $\Psi$ we see from Theorem
\ref{trivthm} and Lemma \ref{curvecount} that
\[N_{\Psi}(2H)\ll H^{2+\ve}.\]
It follows
that the total contribution from all $B$-exceptional lines $L$ with
height $H<H(L)\le 2H$ is
\[\ll \frac{B^2(\log B)^{60}}{H}.H^{2+\ve}.\]
Since we have $H(L)\le B(\log B)^{125}$, by Theorem
\ref{heights-planes}, we may sum this bound over
dyadic intervals with $H\ll B(\log B)^{125}$ to obtain an estimate
$O(B^{3+2\ve})$ for the total contribution to $N(F;B)$ arising from
points on $B$-exceptional lines.  On re-defining $\ve$ this is satisfactory
for Theorem \ref{3fold}.

It now remains to deal with the remaining lines, which we refer to as 
``$B$-ordinary''.  According to our 
definition of $B$-exceptional lines, there will be no $B$-ordinary lines unless
$\dim F_2(X)=1$.  Moreover, when $\dim F_2(X)=1$, any $B$-ordinary line
will be contained in a plane, by Lemma \ref{l-i-p}, but not in a
$B$-good plane.  Moreover the corresponding point
$[\ma{m}_1]$ will not be a vertex of $X$.  Any $B$-ordinary line 
$L=M_j$ with $H<H(L)\le 2H$ will contribute $O(B^2/H)$ to $N(F;B)$, by 
Lemma \ref{planecount}, since 
$$
H(L)\ll B^{1+\ve} 
$$ 
in Theorem \ref{heights-planes}.
Thus we need to estimate the number of $B$-ordinary lines $L$ with height of
order $H$.  To do this we observe that $L$ will be spanned by the
lattice generated by $\ma{m}_1,\ma{m}_2$ with 
\[
H(L)\leq |\ma{m}_1||\ma{m}_2|\ll H(L), 
\]
by (\ref{rog6}) and (\ref{sch}).  In particular we will have  
$|\ma{m}_1|\ll H^{1/2}$, since $s_1\le s_2$. 
We therefore count $B$-ordinary lines according
to the corresponding vector $\ma{m}_1$.  Since $[\ma{m}_1]$ is 
not a vertex of $X$ it follows from Lemma \ref{infplanes} that there 
are $O(1)$ planes through $[\ma{m}_1]$.  If such a plane is not 
$B$-good, it can contain at most one line $M_j$.  To see this we recall
that a line $M_j$ corresponds to a lattice $\sfm_j$ with successive
minima of size $O(B(\log B)^{30})$.  Hence, for large enough $B$, the
line $M_j$ will contain two distinct points of height at most $B(\log
B)^{31}$. However, a plane that is not $B$-good can contain at most two
linearly independent points of such height, and hence can contain 
at most one line $M_j$, as claimed.
We may therefore conclude that there are at most $O(1)$ $B$-ordinary 
lines through each point $[\ma{m}_1]$.  In view of Theorem 
\ref{trivthm}, the number of available $\ma{m}_1$ with 
$|\ma{m}_1|\ll H^{1/2}$ is $O(H^2)$, so that the total
contribution from $B$-ordinary lines $L$ with height $H<H(L)\le 2H$ is
\[\ll \frac{B^2}{H}.H^2=B^2H\ll B^{3+\ve}.\]
On summing over dyadic ranges for $H$, and re-defining $\ve$, we find
that this too is satisfactory for Theorem \ref{3fold}.  This
completes the proof.

\end{document}